\documentclass[12pt]{article}
\usepackage{amsmath,amssymb,amsbsy}
\usepackage{epsfig}
\usepackage{amsfonts}
\usepackage[small]{caption}

\newtheorem{theorem}{Theorem}

\newtheorem{condition}[theorem]{Condition}

\newtheorem{definition}[theorem]{Definition}

\newtheorem{lemma}[theorem]{Lemma}

\newtheorem{problem}[theorem]{Problem}
\newtheorem{proposition}[theorem]{Proposition}

\newenvironment{proof}[1][Proof]{\noindent\textbf{#1.} }{\ \rule{0.5em}{0.5em}}

\renewcommand{\Re}{{\text{Re\,}}}
\renewcommand{\Im}{{\text{Im\,}}}

\newcommand{\RR}{\mathbb{R}}

\newcommand{\ZZ}{\mathbb{Z}}

\newcommand{\zexp}{{\mathcal{Z}_e}}
\newcommand{\zlin}{{\mathcal{Z}_\ell}}
\newcommand{\zprop}{{\mathcal{Z}_p}}
\newcommand{\zpropk}{{\mathcal{Z}_p^\kappa}}
\newcommand{\zpropmk}{{\mathcal{Z}_p^{\text{-}\kappa}}}
\newcommand{\etak}{\eta^\kappa}
\newcommand{\etamk}{\eta^{\text{-}\kappa}}

\newcommand{\eps}{\varepsilon}

\newcommand{\ii}{\mathrm{i}}
\newcommand{\inc}{u^{\mathrm{inc}}}
\newcommand{\uinc}{u^{\mathrm{inc}}}
\newcommand{\ua}{u_\text{ad}}
\newcommand{\uainc}{u_\text{ad}^\text{inc}}
\newcommand{\uascat}{u_\text{ad}^\text{sc}}
\newcommand{\scat}{u^{\mathrm{sc}}}
\newcommand{\uz}{\breve u_0}
\newcommand{\dn}{\partial_n}
\newcommand{\Div}{\nabla\hspace{-2.5pt}\cdot\hspace{-1.5pt}}

\newcommand{\Hone}{{H^1(\Omega)}}

\newcommand{\Honez}{{H^1_0(\Omega)}}

\newcommand{\Honek}{{H^1_{\kappa}(\Omega)}}
\newcommand{\Honemk}{{H^1_{\!\text{-}\kappa}(\Omega)}}

\newcommand{\Ltwo}{{L^2(\Omega)}}

\newcommand{\Tk}{T_\kappa}
\newcommand{\Tmk}{T_{\!\text{-}\kappa}}
\newcommand{\half}{{\textstyle{\frac{1}{2}}}}
\newcommand{\littleo}{\text{\scriptsize ${\cal O}$}}

\newcommand{\trans}{{\cal E}}
\newcommand{\intgp}{\int_{\Gamma_+}\!\!}

\newcommand{\bu}{{\breve u}}
\newcommand{\cu}{{\tilde u}}
\newcommand{\buz}{{\breve u_0}}

\newcommand{\cbuz}{{\tilde{\breve u}_0}}
\newcommand{\btau}{{\breve\tau}}
\newcommand{\beps}{{\breve\eps}}
\newcommand{\ctau}{{\tilde\tau}}
\newcommand{\ceps}{{\tilde\eps}}

\newcommand{\Lnorm}[3]{\| #3 \|_{L^{#1}(#2)}}
\newcommand{\Hnorm}[3]{\| #3 \|_{H^{#1}(#2)}}

\begin{document}

\bibliographystyle{plain}

\title{Field sensitivity to ${L^{}}^p$ variations of a scatterer\footnote{Preprint of an article to be published in J. Math. Anal. Appl., 2009.}}
\author{Stephen P. Shipman\footnote{Dept. of Mathematics, Louisiana State University, Baton Rouge, LA 70803, {\tt shipman@math.lsu.edu}}\\
Louisiana State University}
\maketitle

\begin{abstract}\noindent
For the problem of diffraction of harmonic scalar waves by a lossless periodic slab scatterer, we analyze field sensitivity with respect to the material coefficients of the slab.  The governing equation is the Helmholtz equation, which describes acoustic or electromagnetic fields.  The main theorem establishes the variational (Fr\'echet) derivative of the scattered field measured in the $H^1$ (root-mean-square-gradient) norm as a function of the material coefficients measured in an $L^p$ ($p$-power integral) norm, with $2~<~p~<~\infty$, as long as these coefficients are bounded above and below by~positive constants and do not admit resonance.  The derivative is Lipschitz continuous.  We also establish the variational derivative of the transmitted energy with respect to the material coefficients in $L^p$.
\end{abstract}


\section{Introduction}

\begin{figure}\label{figstructure} 
  \centerline{\scalebox{0.35}{\includegraphics{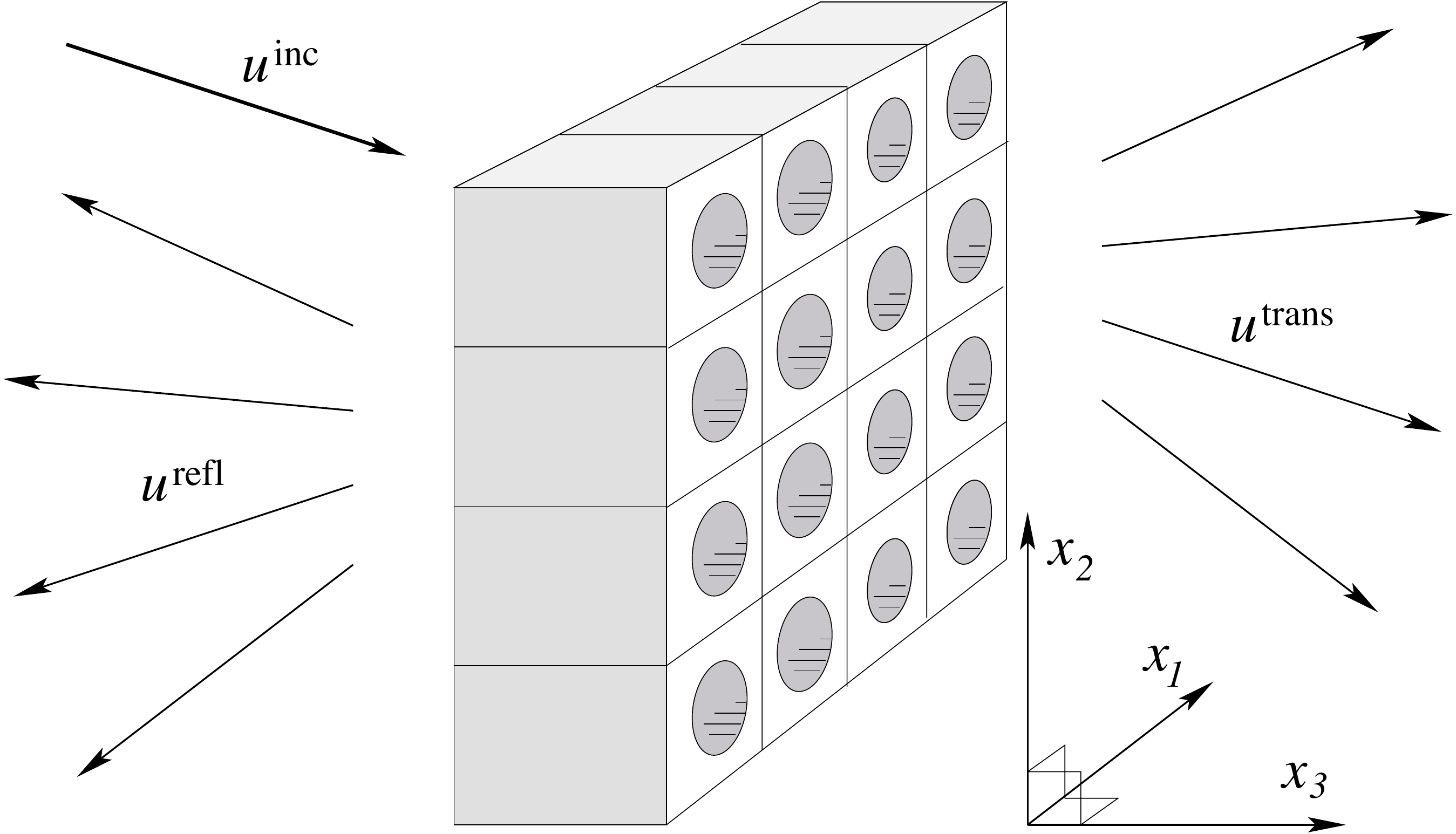}}}
  \caption{An example of a periodic slab scatterer.  The slab is infinite and periodic in the $x_1$ and $x_2$ directions (only sixteen periods are shown) but of finite thickness in the $x_3$ direction.  We impose no fixed boundary conditions (as Dirichlet or Neumann); rather, the slab is in natural contact with the ambient space to its left and right, which is homogeneous with $\eps(x)=\eps_0$ and $\tau(x)=\tau_0$.  This figure depicts a slab consisting of two homogeneous materials.  In this work, we treat the much more general case, in which $\eps(x)$ and $\tau(x)$ are measurable functions bounded from below and above.  The reflected and transmitted fields depict the diffractive orders associated with the angle of incidence of the source field.}
\end{figure}

This work treats the variational calculus of time-harmonic fields scattered by a periodic slab structure as functions of the material coefficients of the scatterer.  We deal with scalar fields~$u$ governed by the linear Helmholtz equation
\begin{equation*}
  \Div \tau \nabla u + \omega^2\eps\, u \,=\, 0,
\end{equation*}
which governs acoustic fields and, in case the coefficients are invariant in one direction, polarized electromagnetic fields, in a composite material characterized by the spatially varying coefficients $\eps$ and $\tau$.  We take these coefficients to be real and positive, which means that the structure is lossless.  Figure \ref{figstructure} depicts an example of the type of scatterer we consider.  The slab is periodic in two directions and finite in the other, and it is in contact with the ambient space, making it an open waveguide.  A~traveling time-harmonic wave, originating from sources exterior to the slab, is incident upon the slab at an angle and is diffracted by~it.  Our aim is to compute the sensitivity of the resulting total field to variations of the material properties ($\eps$ and $\tau$) and geometry of the slab, as well as the sensitivity of the amount of energy transmitted across the slab.

A motivation for this subject is the desire to optimize the way in which energy flows through a periodic slab or film, as well as the related inverse problem, in which one seeks to determine the structure that produces given diffracted field patterns upon illumination by plane waves.  Slabs of photonic crystal structures can be used to guide energy of an incident wave at specific frequencies through channels to the other side of the slab \cite{MoussaWangTuttle2007}.  The characteristics that one seeks to optimize are the amount of transmitted energy and the directionality of the field that is transmitted, as well the electromagnetic mode density, which is important for control of the spontaneous emission rate of atoms placed in the structure \cite{BendicksoDowlingScalora1996}.  The variational calculus of the scattered, or diffracted, field as a function of the structural parameters is the basis for control and optimization of these properties.

Variations of practical interest are not typically uniformly small across the scatterer; rather they tend to be large but supported in a small domain.  For example, one may wish to vary the diameter of a dielectric sphere $S$ of fixed permittivity $\eps_1$ repeated in a two-dimensional periodic array within a matrix of permittivity $\eps_0$ or the diameters of the holes in the example of Figure \ref{figstructure}.  The function $\eps=\eps_1{\chi_{}}_S + \eps_0{\chi_{}}_{S^c}$ (${\chi_{}}_S$ is the characteristic function of $S$) is not continuous with respect to the diameter of $S$ if the function is measured in the supremum norm, or $L^\infty$ norm, but it is continuous if the function is measured in any $p$-power integral norm, or $L^p$ norm, with $1\leq p <\infty$.  Therefore, the question of whether the scattered field is differentiable, or even merely continuous, with respect to $L^p$ perturbations of the scatterer is an important one.  In this work, a rigorous formulation (Theorems \ref{thmmain} and~\ref{thmtrans}) of the following theorem is proved:

\addtocounter{theorem}{-1}
\begin{theorem}\label{thm0}
The scattered field of a lossless periodic slab as well as the transmitted energy, for a fixed incident wave, are Fr\'echet differentiable with respect to the coefficients $\eps$ and $\tau$, if the field and its gradient are measured in the root-mean-square norm (Sobolev norm $H^1$) and $\eps$ and $\tau$ are measured in an ${L^{}}^p$ norm, with $p<\infty$, as long as $\eps$ and $\tau$ are bounded from above and below and do not admit resonance.  The derivatives are Lipschitz continuous.
\end{theorem}

Fr\'echet differentiability with respect to the material coefficients in an $L^p$ norm implies H\"older continuity with respect to variations of a smooth boundary of a homogeneous component of the scatterer.  In fact, the scattered field has been shown to be differentiable with respect to variation of periodic interfaces separating materials of differing dielectric coefficient in two-dimensional polarized electromagnetic scattering problems.  See, for example, Bao \cite{Bao1998}, Dobson \cite{Dobson1999}, and Elschner and Schmidt \cite{ElschnerSchmidt2001,ElschnerSchmidt2003}, as well as \cite{ElschnerSchmidt1998,ElschnerSchmidt2001a}, \cite{Dobson1999}, and Bao and Bonnetier \cite{BaoBonnetier2001} for applications to optimal design.  The differentiability of solutions to {\em strongly} elliptic equations in a bounded domain as well as functionals of these solutions, with respect to the boundary in norms of H\"older continuity, is treated by Pironneau \cite{Pironneau1983} (\S 1.7, 6.2).  In their study of the inverse problem for bounded impenetrable obstacles, Colton and Kress \cite{ColtonKress1998} (\S5.3) prove the Fr\'echet differentiability of the far field pattern in the the $L^2$ norm of the sphere as a function of the boundary in the norm of continuous differentiability.  The inverse problem for scattering by periodic interfaces is treated by Kirsch \cite{Kirsch1994} and in \cite{ElschnerSchmidt2001}.

Theorem~\ref{thm0} implies differentiability with respect to any ${L^{}}^q$ norm with $q\geq p$.  Obtaining an upper bound on the minimal $p$ is an open problem, whose solution would facilitate numerical implementation of the variational gradient.  The formal calculus of variations leads to a candidate for the gradient of the field and transmission coefficient as functions of $\eps$ and~$\tau$.  The gradient of the transmission coefficient is expressed in terms of an adjoint problem, derived formally by Lipton, Shipman, and Venakides \cite{LiptonShipmanVenakides2003}.
In that work, the authors used the formal results in a two-dimensional reduction (where $\eps$~and~$\tau$ are constant in one direction) to manipulate numerically the transmission coefficient as a function of frequency by varying the slab structure.  In this work, this gradient is established rigorously for $L^p$ perturbations of $\eps$ and $\tau$.

The proof uses N. Meyers' theorem on higher integral regularity ($p>2$) of solutions of elliptic equations and their gradients \cite{Meyers1975}.  In order to apply the theorem, one needs an {\it a~priori} bound on the solution of the scattering problem that is independent of the material coefficients.  The obstruction to such a bound is field resonance in the structure, resulting from the presence of guided modes.  A {\em guided mode} is a pseudoperiodic solution (Bloch solution) to the Helmholtz equation that falls off exponentially with distance from the slab.  Mathematically, it is self-sustained, that is, not forced by an incident source field.  Because a solution of the scattering problem is not unique for a given structure at a frequency and Bloch wavevector that admit a guided mode, the scattered field is not uniformly bounded near these parameters.  This work concerns the perturbation of the material properties within a range that excludes resonance, in which the scattering problem necessarily has a unique solution.  Perturbation analysis near resonance is singular, and quite a different problem, as the field and transmitted energy exhibit anomalous behavior near a guided mode frequency.  Rigorous perturbation analysis with respect to frequency and Bloch wavevector about a guided mode is presented in \cite{ShipmanVenakides2005}, and similar analysis with respect to material coefficients and geometry is possible.

\smallskip

The exposition of the ideas and results is summarized as follows.
\begin{itemize}
  \item  Section \ref{secscattering} presents the mathematical formulation of the problem of scattering of incident traveling waves by a periodic slab structure as well as formal perturbation analysis.  This analysis gives rise to the correct candidate for the derivative of the scattered field~$u$ with respect to variations of the material coefficients $\eps$ and $\tau$.
  \item  The sensitivity of the transmitted energy to variations of the scatterer is discussed in Section \ref{sectransmission}, with specialization to structures with homogeneous components.
  \item  Section \ref{seceigenvalues} develops the weak formulation of the scattering problem, in which the frequency and Bloch wavevector are parameters.  We discuss eigenvalues of the sesquilinear forms associated with the scattering problem and their relation to guided modes of the slab.
  \item  The main contributions of this work are stated and proved in Section \ref{secmain}.  The first, Theorem~\ref{thmfieldbound} (p.~\pageref{thmfieldbound}), establishes an {\it a priori} bound on the root-mean-square norm of the solution of the scattering problem and its gradient in the scatterer as long as the structure, frequency, and wavevector do not admit a guided mode.  This result, together with Meyers' regularity theorem are used to prove the main result, Theorem~\ref{thmmain} (p.~\pageref{thmmain}), on the differentiability of the scattered field with respect to $L^p$ variations of the material coefficients.  Theorem~\ref{thmtrans} applies the main theorem and the adjoint method to give an explicit representation of the variational gradient of the transmitted energy as a function of the material coefficients.
\end{itemize}

\section{The scattering problem and sensitivity analysis}\label{secscattering}

The aim of this section is to derive a candidate for the variational gradient of the field scattered by a periodic slab as a function of the material coefficients of the scatterer.  The variational calculus is treated rigorously in Section \ref{secmain}.  First, we present the mathematical formulation of the scattering problem.

\subsection{The scattering problem}

We shall consider time-harmonic solutions $U( x,t) = \tilde u( x)e^{-\ii \omega t}$ ($x=(x_1,x_2,x_3)\in\RR^3$) of the scalar wave equation
\begin{equation}
\eps \frac{\partial^2}{\partial t^2} U = \Div \tau \nabla U,  
\end{equation}
in which the material coefficients $\eps(x)$ and $\tau(x)$ are positive, $2\pi$-periodic in $x_1$ and $x_2$, and bounded from below and above.
The spatial factor $\tilde u$ satisfies the Helmholtz equation
\begin{equation}
\Div \tau \nabla \tilde u + \eps\,\omega^2 \tilde u \,=\, 0.  
\end{equation}
By means of the (partial) Floquet transform in $(x_1,x_2)$, a solution $\tilde u$ can be decomposed into an integral superposition of components $u(x;\kappa)$, where $\kappa=\langle \kappa_1,\kappa_2 \rangle \in\RR^2$,
that are {\em $\kappa$-pseudoperiodic in $x_1$ and~$x_2$ with periods $2\pi$} (see \cite{Kuchment1993} or \cite{ReedSimon1980d}, for example).  This means that $u(x;\kappa)$ satisfies
\begin{eqnarray}
&& \Div \tau \nabla u + \eps\,\omega^2 u \,=\, 0, \label{helmholtz1}\\
&& u( x;\kappa) = u_\text{per}( x;\kappa) e^{\ii(\kappa_1x_1+\kappa_2x_2)} \text{ and $u_\text{per}$ has period $2\pi$ in $x_1$ and $x_2$.}
\end{eqnarray}
The {\em Bloch wavevector} $\kappa=\langle \kappa_1,\kappa_2 \rangle$ is related to the angle of incidence of an incoming wave
\begin{equation}
  e^{i\eta x_3} e^{i((m_1+\kappa_1)x_1 + (m_2+\kappa_2)x_2)},
\end{equation}
for some $m=(m_1,m_2)\in\ZZ^2$, which impinges upon the left side of the slab at an angle of $\theta = \arctan |m + \kappa|/\eta$ with the normal.  We shall take $\kappa$ to lie in the {\em first Brillouin zone},
\[
\kappa \in [-\half,\half)^2,
\]
because each $\kappa\in\RR^2$ can be written $\kappa=m+\bar\kappa$ for $m\in\ZZ^2$ and $\bar\kappa\in[-\half,\half)^2$.

Exterior to the slab ($x_3<z^-$ and $x_3>z^+$), where the material is homogeneous, we set $\eps\!=\!\eps_0>0$ and $\tau\!=\!\tau_0>0$, and the periodic factor $u_\text{per}$ can be decomposed into Fourier components parallel to the slab (in $x'=(x_1,x_2)$).  They are indexed by $ m\!\in\!\ZZ^2$ (with different coefficients on the two sides of the slab), and the $x_3$-dependence of each component is determined by separation of variables in the equation
$\nabla\cdot\nabla u + \omega^2(\eps_0/\tau_0) u = 0$,
\begin{equation}\label{Fourier}
u( x',x_3;\kappa) = \sum_{ m\in\ZZ^2}
\left( c^+_ m \phi_m^+(x_3) + c^-_ m \phi_m^-(x_3) \right) e^{\ii (m+\kappa) x'},
\end{equation}
in which $\phi^\pm_m$ are independent solutions of the ordinary differential equation
$\phi''_m + \eta_m^2\phi_m\!~=~\!0$, where the numbers $\eta_m$ are defined through
\begin{equation}\label{eta1}
  \eta_ m^2 + | m+\kappa|^2 = \omega^2 \eps_0/\tau_0.
\end{equation}
The Fourier harmonics are known as the {\em diffractive orders} or {\em diffraction orders} associated with the periodic structure.  There are many references that expound these ideas, including Wilcox \cite{Wilcox1984} and Nevi\`ere \cite{Neviere}.
The $\phi_m^\pm$ are either oscillatory, linear, or exponential, depending on the numbers $\eta_m$.
We make the following definitions:
\begin{equation}\label{eta}
  \renewcommand{\arraystretch}{1.2}
\left.
  \begin{array}{ccl}
     m\in\zprop \iff \eta_m^2>0, & \eta_m>0 & \text{(propagating harmonics),} \\
     m\in\zlin \iff \eta_m^2=0, & \eta_m=0 & \text{(linear harmonics),} \\
     m\in\zexp \iff \eta_m^2<0, & -i\eta_m>0 & \text{(evanescent harmonics).} \\
  \end{array}
\right.
\end{equation}
%

In the problem of scattering of source fields given by traveling waves impinging upon the slab, we must exclude exponential or linear growth of $u$ as $|x_3|\to\infty$.  The form the of the total field is therefore (we suppress the $\kappa$-dependence of $u(x'\!,x_3;\kappa)$)
\begin{eqnarray}\label{left}
  &u(x'\!,x_3) = \displaystyle\sum_{m\in\zprop} {a_m^\text{inc}e^{i\eta_m x_3}}e^{i(m+\kappa)x'}
  + \sum_{m\in\ZZ^2} a_m e^{-i\eta_mx_3} e^{i(m+\kappa)x'}
  \quad \text{($x_3\leq z_-$),}& \\ \label{right}
  &u(x'\!,x_3) = \displaystyle\sum_{m\in\zprop} {b_m^\text{inc}e^{-i\eta_m x_3}}e^{i(m+\kappa)x'}
  + \sum_{m\in\ZZ^2} b_m e^{i\eta_mx_3} e^{i(m+\kappa)x'}
  \quad \text{($x_3\geq z_+$).}&
\end{eqnarray}
The infinite series are understood in the $L^2$ sense.  The first sums in these expressions represent right-traveling source waves incident upon the slab from the left side and left-traveling source waves incident upon the slab from the right side.  We say that a function $u$ is {\em outgoing} if it is of the form (\ref{left},\ref{right}) with $a_m^\text{inc}=0$ and $b_m^\text{inc}=0$ for all $m\in\zprop$.

\begin{definition}[Outgoing and incoming]
A complex-valued function $u$ defined on $\RR^3$ is said to be {\em outgoing} if there are real numbers $z_-$ and  $z_+$ and sequences $\{a_m\}_{-\infty}^\infty$ and $\{b_m\}_{-\infty}^\infty$ in $\ell^2(-\infty,\infty)$ such that
\begin{eqnarray}  \label{left1}
  & u(x) =  \displaystyle\sum_{m\in\ZZ^2} a_m e^{-i\eta_mx_3} e^{i(m+\kappa)x'}
  \quad \text{for  $x_3\leq z_-$,}& \\ \label{right1}
  &u(x) = \displaystyle\sum_{m\in\ZZ^2} b_m e^{i\eta_mx_3} e^{i(m+\kappa)x'}
  \quad \text{for $x_3\geq z_+$.}&
\end{eqnarray}
The function $u$ is said to be {\em incoming} if it admits the expansions
\begin{eqnarray}  \label{left1in}
  & u(x) =  \displaystyle\sum_{m\in\ZZ^2} a_m e^{i\eta_mx_3} e^{i(m+\kappa)x'}
  \quad \text{for  $x_3\leq z_-$,}& \\ \label{right1in}
  &u(x) = \displaystyle\sum_{m\in\ZZ^2} b_m e^{-i\eta_mx_3} e^{i(m+\kappa)x'}
  \quad \text{for $x_3\geq z_+$.}&
\end{eqnarray}
\end{definition}

We shall take the pseudoperiodic source field to be a superposition of traveling waves incident upon the slab from the left and right.  We think of these waves as emanating from $x_3=-\infty$ and from $x_3=\infty$:
\begin{equation}\label{incident}
\inc(x'\!,x_3) = \sum_{m\in\zprop} (a_m^\text{inc}e^{i\eta_m x_3} + b_m^\text{inc}e^{-i\eta_m x_3})
e^{i(m+\kappa)x'}.
\end{equation}
The problem of scattering of the incident wave $\inc$ by the slab is expressed as a system characterizing the {\em total field} $u$, which is the sum of the incident field $\uinc$ and the {\em scattered, or diffracted, field} $\scat$, the latter of which is {outgoing}.  The ``strong form" of the problem is posed for functions $\eps$ and $\tau$ that are smooth except on a set $\Sigma$ consisting of continuously differentiable surfaces of discontinuity, with normal vector $n$.  The ``weak form", presented in subsection \ref{subsecweak}, allows $\eps$ and $\tau$ to be merely measurable.
\begin{problem}[Scattering of an incident wave, strong form] 
\label{problemscattering1}
Given an incident field \eqref{incident}, find a function $u$ that satisfies the following conditions.
\begin{eqnarray}
&& \Div \tau \nabla  u + \omega^2\eps  u \,=\, 0 \;\; \text{in $\RR^3\setminus\Sigma$}, \label{scattering1}\\
&& u \text{ and } \tau \dn u \text{ are continuous on $\Sigma$,}\label{scattering2}\\
&& \text{$u$ is $\kappa$-pseudoperiodic in $(x_1,x_2)$},\label{scattering3}\\
&& u = \inc + \scat, \text{ with $\scat$ outgoing.}\label{scattering4}
\end{eqnarray}
\end{problem}

One can generalize the scattering problem by introducing sources originating from the slab itself or from points outside the slab.  Such sources are represented by a periodic function $h$ and a periodic vector field $\xi$, which enter the equation thus:
\begin{equation}
  \Div \tau \nabla  u + \omega^2\eps  u \,=\, \Div\xi + h.
\end{equation}
In our investigation of the perturbation of the scattering Problem \ref{problemscattering1}, we will be concerned with an auxiliary problem involving sources that are confined to the region between $x_3=z_-$ and $x_3=z_+$ (besides the incident source field originating from $x_3=\pm\infty$).
\begin{problem}[General scattering, strong form] 
\label{problemscatteringgen1}
Given an incident field \eqref{incident}, find a function $u$ that satisfies the following conditions.
\begin{eqnarray}
&& \Div \tau \nabla  u + \omega^2\eps  u \,=\, \Div\xi + h \;\; \text{for } z_-<x_3<z_+ \text{ and } x\notin\Sigma, \label{scatteringgen0}\\
&& \Div \tau \nabla  u + \omega^2\eps  u \,=\, 0 \;\; \text{otherwise}, \label{scatteringgen1}\\
&& u \text{ and } \tau \dn u \text{ are continuous on $\Sigma$,}\label{scatteringgen2}\\
&& \text{$u$ is $\kappa$-pseudoperiodic in $(x_1,x_2)$},\label{scatteringgen3}\\
&& u = \inc + \scat, \text{ with $\scat$ outgoing.}\label{scatteringgen4}
\end{eqnarray}
\end{problem}

\subsection{Formal sensitivity analysis}\label{secformal}

Let $u$ be the solution of the scattering Problem \ref{problemscattering1} (existence and uniqueness will be dealt with later), and let $u+\breve u$ be the solution of the scattering problem with the same incident field but with $\eps+\breve\eps$ and $\tau+\breve\tau$ in place of $\eps$ and~$\tau$.  The coefficients $\eps_0$ and $\tau_0$ exterior to the slab remain fixed.  The functions $u$ and $u+\bu$ satisfy
\begin{eqnarray}
  & \Div\tau\nabla u + \omega^2\eps u = 0, & \\
  & \Div(\tau+\breve\tau)\nabla(u+\breve u) + \omega^2(\eps+\breve\eps)(u+\breve u) = 0. &
\end{eqnarray}
Subtracting these equations yields the equation for the perturbed field $\breve u$,
\begin{equation}\label{PDEperturbed}
  \Div(\tau+\breve\tau)\nabla\breve u + \omega^2(\eps+\breve\eps) \breve u
    = -\Div\breve\tau\nabla u + \omega^2\breve\eps u,
\end{equation}
and $\breve u$ is outgoing because the incident fields for $u$ and $u+\breve u$ are identical and the forcing term on the right-hand side of \eqref{PDEperturbed} is confined to $z_-\!<\!x_3\!<\!z_+$.  If we remove the terms on the left-hand side that are quadratic in $\breve\eps$, $\breve\tau$, and $\breve u$, we obtain the differential equation for the formal leading-order sensitivity $\breve u$ of the total field as a function of the perturbations $\breve\eps$~and~$\breve\tau$.  We denote this linear approximation to $\breve u$ by $\breve u_0$,
\begin{eqnarray}
  && \Div\tau\nabla\breve u_0 + \omega^2\eps \breve u_0
  = -\Div\breve\tau\nabla u + \omega^2\breve\eps u, \\
  && \breve u_0 \text{ is outgoing.}
\end{eqnarray}

In order to establish that the linear map $(\breve\eps,\breve\tau)\mapsto\breve u_0$ is truly the variational differential of $u$ with respect to $(\eps,\tau)$, we should demonstrate two things,
\begin{eqnarray}
  & \| \breve u_0 \| \leq C \|(\breve\eps,\breve\tau)\|, & \label{differential1} \\
  & \displaystyle\frac{\| \breve u - \breve u_0 \|}{\|(\breve\eps,\breve\tau)\|} \to 0 \text{ as } \|(\breve\eps,\breve\tau)\|\to 0.& \label{differential2}
\end{eqnarray}
The appropriate norm in which to measure $u$ is the $H^1$-norm, restricted to a domain $\Omega$ (Figure 2) comprising one period of the structure between $x_3=z_-$ and $x_3=z_+$:
\begin{eqnarray}
  && \Omega = \{ x\in\RR^3 : 0<x_1<2\pi,\, 0<x_2<2\pi,\, z_-<x_3<z_+ \}, \label{Omega}\\
  && \Gamma_\pm = \{ x\in\RR^3 : 0<x_1<2\pi,\, 0<x_2<2\pi,\, x_3=z_\pm \}. \label{Gamma}
\end{eqnarray}
The two-dimensional squares $\Gamma_\pm$ are the left and right boundaries of $\Omega$.
The normal vector $n$ to $\Gamma$ is taken to be directed outward, so that
\begin{equation}
  \partial_n u =
  \renewcommand{\arraystretch}{1.1}
\left\{
  \begin{array}{ll}
	-{\partial u}/{\partial x_3} & \text{on } \Gamma_-, \\
	{\partial u}/{\partial x_3} & \text{ on } \Gamma_+.
  \end{array}
\right.
\end{equation}
The $H^1$ norm in $\Omega$~is
\begin{equation}
  \Hnorm{1}{\Omega}{u} = \left( \int_\Omega (|\nabla u|^2 + |u|^2)dV \right)^{1/2}.
\end{equation}
The main result of this work proves that \eqref{differential1} and \eqref{differential2} hold if $(\breve\eps,\breve\tau)$ is measured in some $L^p$ norm in $\Omega$, with $p<\infty$,
\begin{equation}
  \|(\eps,\tau)\|_{L^p(\Omega)} = \left( \int_\Omega (|\breve\eps|^p + |\breve\tau|^p) dV \right)^{1/p}.
\end{equation}

\begin{figure}\label{figdomain} 
  \centerline{\scalebox{0.4}{\includegraphics{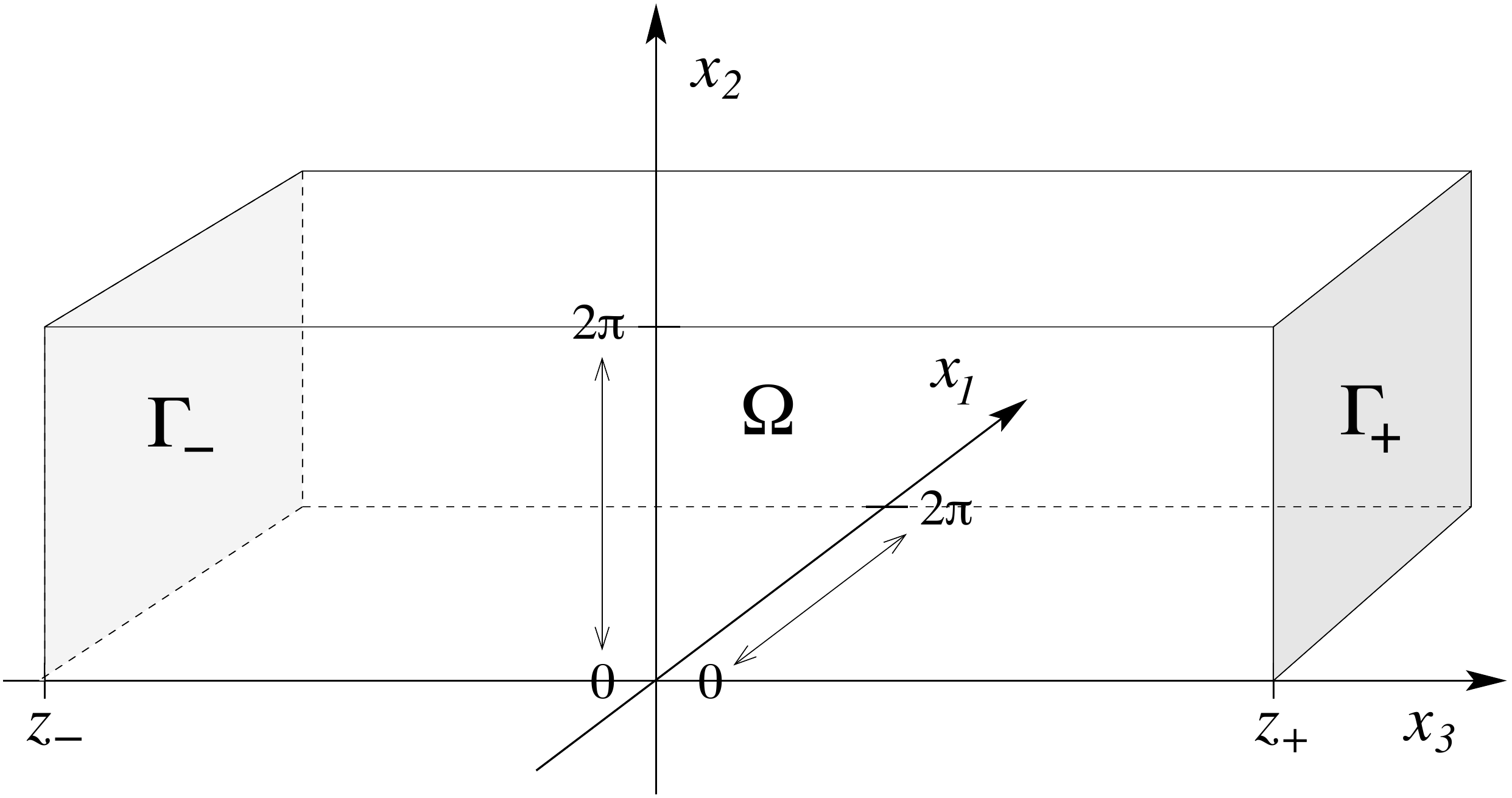}}}
  \caption{The rectangular prism $\Omega$ encloses one period of the slab structure, which is defined by material coefficients $\eps(x)$ and $\tau(x)$ for $x=(x_1,x_2,x_3)\in\Omega$ that are extended periodically in the $x_1$ and $x_2$ directions.  For $x_3<z_-$ and for $x_3>z_+$, the material is homogeneous, with $\eps(x)=\eps_0$ and $\tau(x)=\tau_0$.}
\end{figure}

\section{Energy transmission and special structures}\label{sectransmission}

We apply the results discussed in the previous section to the sensitivity analysis of the amount of energy of an incident wave that is transmitted from one side of the slab to the other.  The formal analysis of the adjoint problem associated with the differential of the transmitted energy that was derived in \cite{LiptonShipmanVenakides2003} is revisited in the light of the rigorous results of this work.

\subsection{Variation of the transmitted energy}

Let us send a traveling wave toward the slab from the left and consider the energy transmitted to right side of the slab.  This means that we take $b^\text{inc}_m=0$ in \eqref{right}.  We are interested in the sensitivity of the transmitted energy to perturbations of the material coefficients $\eps$ and $\tau$.  The time-averaged energy flux through one period of the right-hand boundary of the slab is defined by
\begin{equation}
  \trans = \Im\intgp \tau_0\, \bar u \partial_n u,
\end{equation}
in which $u$ is the solution to the scattering Problem \ref{problemscattering1}.  This quantity can be expressed in terms of the Fourier coefficients $b_m$ of the propagating harmonics of the transmitted field:
\begin{eqnarray}
  && \displaystyle u = \sum_{m\in\ZZ^2} b_m e^{i\eta_mx_3} e^{i(m+\kappa)x'}
  \quad \text{($x_3\geq z_+$),} \\
  && \trans = \tau_0 \sum_{m\in\zprop}\eta_m|b_m|^2.
\end{eqnarray}
The coefficients $b_m$ are functions of $\eps$ and $\tau$.

We prove in section \ref{subsectrans} (Theorem~\ref{thmtrans}) that $\trans$ is differentiable with respect to $\eps$ and $\tau$ if these are measured in an appropriate $L^p$ norm, with $p<\infty$, as long as $\eps$ and $\tau$ are bounded from below and above by positive numbers and there are no resonant frequencies for the scattering problem.  The derivative is expressed in terms of the solution $\ua$ to an adjoint problem in which the incident field $\uainc$ is obtained by sending the transmitted field of $u$ back toward the slab from the right.  The incident and scattered fields have Bloch wavevector $-\kappa$.
\begin{eqnarray}\label{adjointproblem}
&& \Div \tau \nabla  \ua + \omega^2\eps  \ua \,=\, 0 \;\; \text{in $\RR^3\setminus\Sigma$}, \label{ascattering1}\\
&& \ua \text{ and } \tau \dn \ua \text{ are continuous on $\Sigma$,}\label{ascattering2}\\
&& \text{$\ua$ is $-\kappa$-pseudoperiodic in $(x_1,x_2)$,}\label{ascattering3}\\
&& \ua = \uainc + \uascat, \text{ with $\uascat$ outgoing,}\label{ascattering4}\\
&& \displaystyle \uainc = \sum_{m\in\zprop} \bar b_{m} e^{-i\eta_mx_3}e^{-i(m+\kappa)x'}. \label{ascattering5}
\end{eqnarray}
%
If $u+\breve u$ is the solution of the scattering problem with perturbed coefficients $\eps+\breve\eps$ and $\tau+\breve\tau$ and $\breve\trans$ is the corresponding change in the transmitted energy, that is,
\begin{equation}
  \trans + \breve\trans
  = \Im\intgp \tau_0\, (\bar u+\bar{\breve u}) \partial_n (u+\breve u),
\end{equation}
then the linear, leading-order, change $\breve\trans_0$ in $\trans$ is given by 
\begin{equation}\label{trans00}
  \breve\trans_0 = \Im \int_\Omega \left( \breve \tau \, \nabla u \cdot \nabla \ua - \omega^2\, \breve\eps\, u \ua  \right).
\end{equation}
It is proved in Theorem~\ref{thmtrans} that $\breve\trans_0$ is bounded as a function of $\breve\eps$ and $\breve\tau$ measured in an $L^p$ norm ($p<\infty$) and that the error is estimated by the square of the $L^p$ norm,
\begin{equation}\label{trans00est}
  |\breve\trans - \breve\trans_0|
  \leq C \left( \Lnorm{p}{\Omega}{\breve\tau} + \Lnorm{p}{\Omega}{\breve\eps} \right)^2.
\end{equation}

\subsection{Variation of the complex transmission coefficients}

The transmitted energy $\trans$ is a function of the coefficients of the transmitted propagating harmonics $b_m$, $m\in\zprop$.  In fact, one can obtain the variational gradient of each complex coefficient individually.  The associated adjoint problem for $b_m$ is obtained by replacing the incident field \eqref{ascattering5} by a single incoming harmonic,
\begin{equation}
  (\ua^m)^\text{inc} = e^{-i\eta_mx_3}e^{-i(m+\kappa)x'}.
\end{equation}
The linear, leading-order, change $(\breve b_m)_0$ in $b_m$ as a function of variations of $\eps$ and $\tau$ is
\begin{equation}\label{b0}
  (\breve b_m)_0 = \frac{i}{8\pi^2\eta_m\tau_0} \int_\Omega \left( \breve \tau \, \nabla u \cdot \nabla \ua^m - \omega^2\, \breve\eps\, u \ua^m  \right).
\end{equation}
Compare the formulas in \cite{ElschnerSchmidt1998} and \cite{ElschnerSchmidt2003}, \S4, for the case of conical diffraction by two-dimensional periodic structures, in which the interfaces between contrasting homogeneous dielectrics are varied.

\subsection{Structures with homogeneous components}

An important class of periodic structures is comprised of those that consist of homogeneous components.  
The variational gradient \eqref{trans00} can be formulated in terms of the material and geometric parameters of these components.  Suppose that one period of the slab consists of components described by $N$ disjoint domains $D_j$ with material coefficients given by spatial constants $\eps_j$ and $\tau_j$.  The coefficients exterior to the components are $\eps_0$ and $\tau_0$, and the normal vector $n$ to $\partial D_j$ is directed outward.

\smallskip
{\sl Variation of the values of $\tau_j$ and $\eps_j$.} \ If we keep the boundaries of the domains $D_j$ fixed and perturb the numbers $\eps_j$ and $\tau_j$ by amounts $\breve\eps_j$ and $\breve\tau_j$, then \eqref{trans00} becomes
\begin{multline}\label{variationvalues}
  \breve\trans_0 = \Im \sum_{j=1}^N
  \left( \breve\tau_j \int_{D_j} \!\nabla u\cdot\nabla\ua \,-\, \omega^2\breve\eps_j \int_{D_j}\!\!u\ua
  \right) \\
  = \Im \sum_{j=1}^N
  \left( \breve\tau_j \int_{\partial D_j} \!\!\ua\partial_nu
  \,\,+\,\, \omega^2(\breve\tau_j\eps_j\tau_j^{-1}-\breve\eps_j) \int_{D_j}\!\!u\ua
  \right).
\end{multline}
Since the domains $D_j$ are fixed, the estimate \eqref{trans00est} yields
\begin{equation}
  |\breve\trans-\breve\trans_0| \leq
  C \sum_{j=1}^N (|\breve\eps_j|^2+|\breve\tau_j|^2),
\end{equation}
and therefore \eqref{variationvalues} gives the gradient of $\trans$ with respect to the numbers $\eps_j$ and $\tau_j$.

\smallskip
{\sl Variation of the boundaries.} \ Let us now hold $\eps_j$ and $\tau_j$ fixed and let each boundary $\partial D_j$ vary in the direction of a given vector field $v_j$ defined on $\partial D_j$ by allowing the points on $\partial D_j$ to flow in the direction of $v_j$ for a distance $h$.  Then \eqref{trans00} becomes
\begin{equation}
  \breve\trans_0 = \Im\sum_{j=1}^N \int_{\breve D_j}
  \left( \pm(\tau_j-\tau_0) \nabla u\cdot\nabla\ua \mp \omega^2(\eps_j-\eps_0)u\ua
  \right),
\end{equation}
in which $\breve D_j$ denotes the region traversed by the boundary points.  The upper sign is taken if the normal component of $v_j$ is directed out of $D_j$ and the lower sign is taken if the normal component of $v_j$ is directed into $D_j$.  For small $h$, we obtain
\begin{equation}\label{variationboundaries}
  \breve\trans_0
  \,=\, h \,\Im \sum_{j=1}^N \int_{\partial D_j}
  \left( (\tau_j-\tau_0)\nabla u\cdot\nabla\ua - \omega^2(\eps_j-\eps_0) u\ua
  \right) v_j\!\cdot\! n \,+\, \littleo(h).
\end{equation}
From \eqref{trans00est}, we obtain, for sufficiently small $h$,
\begin{equation}
  |\breve\trans-\breve\trans_0| \leq
  C\, h^{2/p} \left( \sum_{j=1}^N
    \int_{\partial D_j} (|\eps_j-\eps_0|^p + |\tau_j-\tau_0|^p) |v_j\!\cdot\!n|
  \right)^{2/p}.
\end{equation}
This result implies only that $\trans$ is differentiable with respect to $h^r$ for $0<r<2/p$ at $h=0$; in particular, $\trans$ is H\"older continuous with respect to uniform perturbations of the boundary.  As discussed in the Introduction, it has been proven in two-dimensional cases that $\trans$ is in fact differentiable with respect to $h$.

\section{Eigenvalues and the scattering problem in weak form}\label{seceigenvalues} 

The weak formulation of the scattering problem places it within the framework of sesquilinear forms in the Hilbert space $\Honek$ of $\kappa$-pseudoperiodic functions on a period $\Omega$ of the scatterer.  It allows proper treatment of guided modes, as well as existence, uniqueness, and bounds of solutions.  The weak formulation requires the Dirichlet-to-Neumann map that characterizes outgoing fields.  For bounded scatterers in $\RR^3$, one may refer to Lenoir, {\em et. al.},~\cite{LenoirVulliermeHazard1992} or Colton and Kress \cite{ColtonKress1998} (\S 5.3); for periodic structures, our formulation is essentially the same as that used by Bonnet-Bendhia and Starling \cite{Bonnet-BeStarling1994}.  

\subsection{The weak formulation of the scattering problem}\label{subsecweak}

By treating the Helmholtz equation in the scattering Problem \ref{problemscattering1} in the weak sense, the second condition on the the continuity of $u$ and $\tau\partial_n u$ is automatically satisfied. The weak sense is expressed as follows: If $\eps$ and $\tau$ are smooth except along smooth surfaces of discontinuity, and if $\xi$ is a smooth vector field and $h$ is a smooth scalar function, then a function $u$ satisfies the Helmholtz equation
\begin{equation}
  \Div \tau \nabla  u + \eps\,\omega^2  u \,=\, \Div\xi + h
\end{equation}
at points where $\eps$ and $\tau$ are smooth {\em and} the condition of continuity of $u$ and $\tau\partial_n u$ on interfaces between materials if and only if
\begin{equation}\label{weaksense}
\int_{\RR^3}
\left( \tau \, \nabla u \cdot \nabla \bar v - \omega^2\, \eps\, u \bar v  \right)
\,=\, 
\int_{\RR^3}
\left( \xi \cdot \nabla \bar v + h \bar v  \right)
  \quad \text{for all } v\in C_0^\infty(\RR^3).
\end{equation}
This weak form of the Helmholtz equation allows one to relax the regularity of $\eps$ and $\tau$ so that they are merely measurable and the regularity of $\xi$, $h$, $u$, and the distributional gradient of $u$, so that they are required only to be locally square-integrable.

To incorporate the pseudoperiodicity and outgoing conditions required by the scattering Problem \ref{problemscattering1}, its weak form is posed in one period $\Omega$ of the slab structure, between its bounding planes $x_3=z_-$ and $x_3=z_+$ (see Fig. 2).
The pseudoperiodicity condition is enforced by requiring that the solution $u$ and the test functions $v$ be in the pseudoperiodic Sobolev space
\begin{equation}
  \Honek = \{ u\in\Hone : u(2\pi,x_2,x_3) = e^{2\pi i\kappa_1}u(0,x_2,x_3), u(x_1,2\pi,x_3) = e^{2\pi i\kappa_2}u(x_1,0,x_3) \}.
\end{equation}
The evaluation of $u$ on the boundary of $\Omega$ is in the sense of the trace map $\Hone\to H^{\half}(\partial\Omega)$.

The outgoing condition is enforced through the Dirichlet-to-Neumann operator for outgoing fields, $T:H^{\half}(\Gamma) \to H^{-\half}(\Gamma)$.
It acts on traces on $\Gamma$ of functions in $\Honek$ and is defined through the Fourier transform as follows.
For any $f\in H^{\half}(\Gamma)$, let $\hat f_ m$ be the Fourier coefficients of $e^{-i\kappa x'}f$; this is a pair of numbers $\hat f_m = (\hat f_m^-,\hat f_m^+)$, one giving the $m^\text{th}$ pseudoperiodic Fourier component of $f$ on $\Gamma_-$ and the other on $\Gamma_+$,
\begin{equation}\label{Fourierbdy}
  f(x_1,x_2,z_\pm) = \sum_{m\in\ZZ^2} \hat f_m^\pm e^{i(m+\kappa)x'}.
  \end{equation}
Then $T$ is defined by
\begin{equation}\label{T}
  T:H^{\half}(\Gamma) \to H^{-\half}(\Gamma), \quad  (\widehat{Tf})_ m = -i\eta_ m \hat f_ m.
\end{equation}
The operator $T$ has a nonnegative real part $T_r$ and a nonpositive imaginary part $T_i$:
\begin{eqnarray}\label{TrTi}
  && T = T_r + i T_i, \\
  && (\widehat{T_r f})_m = \renewcommand{\arraystretch}{1}
\left\{
  \begin{array}{ll}
	-i\eta_m \hat f_m & \text{if } m\in\zexp, \\
	0 & \text{otherwise}.
  \end{array}
\right.  \label{Tr}\\
   && (\widehat{T_i f})_m = \renewcommand{\arraystretch}{1}
\left\{
  \begin{array}{ll}
	-\eta_m \hat f_m & \text{if } m\in\zprop, \\
	0 & \text{otherwise}.
  \end{array}
\right. \label{Ti}
\end{eqnarray}
The adjoint of $T$ with respect to the pairing $(f,g) = \int_\Gamma f\bar g$, for $f\in H^{\half}(\Gamma)$ and $g\in H^{-\half}(\Gamma)$~is
\begin{equation}
  T^*:H^{\half}(\Gamma) \to H^{-\half}(\Gamma), \quad  (\widehat{T^*f})_ m
  = i\bar\eta_ m \hat f_ m.
\end{equation}
$T$ characterizes the normal derivative of an outgoing function on $\Gamma$ as a function of its values on $\Gamma$.  If we denote the trace of $u$ on $\Gamma$ by $u$ again, then
\begin{equation}\label{outgoing}
  \dn u + T u = 0  \text{ on $\Gamma$} \quad \text{for $u$ outgoing},
\end{equation}
whereas the adjoint $T^*$ characterizes incoming fields,
\begin{equation}\label{incoming}
  \dn u + T^* u = 0  \text{ on $\Gamma$} \quad \text{for $u$ incoming}.
\end{equation}
Using this together with the decomposition $u=\inc+\scat$ of the solution to the scattering Problem \ref{problemscattering1}, we obtain
\begin{equation}\label{forcing1}
  \partial_n u + T u \,=\, \partial_n\inc + T\inc \,=\, 
  \renewcommand{\arraystretch}{1.4}
\left\{
  \begin{array}{ll} \displaystyle
	\sum_{m\in\zprop} -2i\eta_m a^\text{inc}_m e^{i\eta_m x_3}e^{i(m+\kappa)x'}, &
	  x\in\Gamma_-, \\ \displaystyle
	\sum_{m\in\zprop} -2i\eta_m b^\text{inc}_m e^{-i\eta_m x_3}e^{-i(m+\kappa)x'}, &
	  x\in\Gamma_+.	
  \end{array}
\right.
\end{equation}
The function $(\partial_n+T)\inc$ gives rise to an element of the space $\Honek^*$ of bounded conjugate-linear functionals on $\Honek$, which we denote by $f^\omega_{\Gamma}$.  We emphasize only the dependence on the frequency $\omega$, the parameters $\kappa$, $\eps_0$, and $\tau_0$ being fixed.  We also write $T^\omega$ for $T$.
\begin{equation}
  f^\omega_{\Gamma}(v) = \tau_0\int_\Gamma ((\partial_n + T^\omega)\inc)\bar v\,dA \quad\text{for all } v\in\Honek,
\end{equation}
in which evaluation of $\bar v$ on $\Gamma$ is in the sense of the trace map.

\begin{problem}[Scattering of an incident wave, weak form]
\label{problemscattering2}
Find a function $u\in\Honek$ such that
\begin{equation}\label{weak}
\int_\Omega
\left( \tau \nabla u \cdot \nabla \bar v - \omega^2 \eps u \bar v  \right)
+ \tau_0 \int_\Gamma (T^\omega u) \bar v
\,=\, f^\omega_{\Gamma}(v)
\quad \text{for all $v\in\Honek$.}  
\end{equation}
\end{problem}
The scattering problem is generalized by allowing $f^\omega_{\Gamma}$ to be replaced by a general element $f\in\Honek^*$.  Problem \ref{problemscatteringgen1} has the weak form
\begin{problem}[General scattering, weak form]
\label{problemscatteringgen2}
Find a function $u\in\Honek$ such that
\begin{equation}\label{weakgeneral}
\int_\Omega
\left( \tau \nabla u \cdot \nabla \bar v - \omega^2 \eps u \bar v  \right)
+ \tau_0 \int_\Gamma (T^\omega u) \bar v
\,=\, f^\omega_{\Gamma}(v) + \int_\Omega \left( \xi \cdot \nabla \bar v + h \bar v  \right)
\quad \text{for all $v\in\Honek$.}  
\end{equation}
\end{problem}
The vector field $\xi$ and the function $h$ are in $\Ltwo$, making the right-hand side a bounded conjugate-linear functional on $\Honek$.
%
%
%
%
%
The equivalence of the scattering Problems \ref{problemscattering1} and \ref{problemscattering2} as well as their generalizations \ref{problemscatteringgen1} and \ref{problemscatteringgen2} is expressed in the following theorem, whose proof is standard.

\begin{proposition}[Equivalence of strong and weak forms]\label{propequivalence}
Let $\eps$ and $\tau$ be bounded and measurable in $\Omega$, and let $\xi$ and $h$ be in $\Ltwo$.
If $u\in H^1_\text{loc}(\RR^3)$ satisfies the scattering Problem~\ref{problemscattering1} (resp. \ref{problemscatteringgen1}), in which the Helmholtz equation {\rm (\ref{scatteringgen0},\ref{scatteringgen1})} and the interface conditions \eqref{scatteringgen2} together are replaced by the weak condition \eqref{weaksense}, then $u|_\Omega \in \Honek$ satisfies Problem \ref{problemscattering2} (resp.~\ref{problemscatteringgen2}).
  Conversely, if $u\in\Honek$ satisfies Problem \ref{problemscattering2} (resp. \ref{problemscatteringgen2}), then there exists a unique extension $\tilde u$ of $u$ to $\RR^3$ such that $\tilde u$ satisfies Problem \ref{problemscattering1} (resp. \ref{problemscatteringgen1}).
\end{proposition}
The unique extension $\tilde u$ of the solution $u$ of Problem \ref{problemscattering2} or \ref{problemscatteringgen2} to all of space, mentioned in Proposition~\ref{propequivalence}, admits the Fourier expansions (\ref{left},\ref{right}).  Because of this, one can prove that $\tilde u$ is bounded in any finite domain in $\RR^3$ by $u$ and the incident field, as expressed in the following theorem.  The theorem may be proved most elegantly using an integral representation formula that expresses the scattered field in a finite region to the left or right of one period of the structure as a bounded operator of the Cauchy data on $\Gamma$ of the total field (see Lemma 2.1 of \cite{RamdaniShipman2008} and the proof of Lemma 3.8 of \cite{CostabelStephan1985} for boundedness).  We take a more direct approach here.

\begin{lemma}[Boundedness of field extension]\label{lemmaextension}
  Let $D$ be a bounded domain in $\RR^3$ and $\omega_+$ a positive number.  There exist numbers $C_1$ and $C_2$, independent of $\omega$ as long as $\omega\leq\omega_+$, such that, if $u\in H^1_\text{loc}(\RR^3)$ is $\kappa$-pseudoperiodic in $x'=(x_1,x_2)$ and admits expansions (\ref{left},\ref{right}), then
\begin{equation}
  \Hnorm{1}{D}{u} \leq C_1 \Hnorm{1}{\Omega}{u} + C_2\!\!\sum_{m\in\zprop} (|a_m^\text{\rm inc}|^2 + |b_m^\text{\rm inc}|^2)(1+|m|).
\end{equation}
\end{lemma}

\begin{proof}
Because $u$ is pseudoperiodic, it is sufficient to prove the theorem for a domain $\tilde\Omega$ of the form
\begin{equation}
  \tilde\Omega = \{ x\in\RR^3 : 0<x_1<2\pi,\, 0<x_2<2\pi,\, z_0<x_3<z_- \},
\end{equation}
and for a domain of an analogous form with $z_+<x_3<z_0$.  The proofs are analogous.
It is convenient to express the form \eqref{left} as
\begin{equation}
    u(x'\!,x_3) = \displaystyle\sum_{m\in\zprop} {a_m^\text{inc}e^{i\eta_m (x_3-z_-)}}e^{i(m+\kappa)x'}
  + \sum_{m\in\ZZ^2} a_m e^{-i\eta_m(x_3-z_-)} e^{i(m+\kappa)x'}
  \quad \text{($x_3\leq z_-$).}
\end{equation}
Denote the first sum by $u_1$ and the second by $u_2$.
\begin{multline}\label{estu2}
  \int_{\tilde\Omega} |u_2|^2 = 4\pi^2 \sum_{m\in\ZZ^2} \int_{z_0}^{z_-} | a_m e^{-i\eta_m(x_3-z_-)} |^2 dx_3
  = 4\pi^2 \sum_{m\in\ZZ^2} |a_m|^2 \int_{z_0}^{z_-} |e^{-2i\eta_m(x_3-z_-)}|dx_3 \\
  \leq 4\pi^2(z_--z_0)\! \sum_{m\not\in\zexp} |a_m|^2 + 4\pi^2\! \sum_{m\in\zexp} \frac{|a_m|^2}{2|\eta_m|}.
\end{multline}
The gradient of $u_2$ in $\tilde\Omega$ is
\begin{equation}
  \nabla u_2 = \sum_{m\in\ZZ^2} a_m e^{-i\eta_m(x_3-z_-)} e^{i(k+m)x'}
  \langle i(\kappa+m),-i\eta_m \rangle. 
\end{equation}
Similar estimates yield
\begin{equation}\label{estgradu2}
  \int_{\tilde\Omega} |\nabla u_2|^2 \leq 4\pi^2(z_--z_0)\! \sum_{m\not\in\zexp} |a_m|^2(\eta_m^2+|\kappa+m|^2)
  \,+\, 4\pi^2\! \sum_{m\in\zexp} |a_m|^2 \frac{|\eta_m|^2+|\kappa+m|^2}{2|\eta_m|}.
\end{equation}
For $u_1$, we obtain the estimate
\begin{equation}\label{estu1}
  \int_{\tilde\Omega} (|u_1|^2 + |\nabla u_1|^2) \leq
  4\pi^2(z_--z_0)\! \sum_{m\in\zprop} |a_m^\text{inc}|^2 (1+\eta_m^2 + (\kappa+m)^2).
\end{equation}
From the estimates (\ref{estgradu2},\ref{estu1}) and the definition of the numbers $\eta_m$, one infers that there is a positive constant $c$ such that
\begin{equation}\label{estu}
 c \int_{\tilde\Omega} (|u|^2 + |\nabla u|^2) \leq \sum_{m\in\zprop} |a^\text{inc}_m|^2(1+|m|)+ \sum_{m\in\ZZ^2} |a_m|^2(1+|m|),
\end{equation}
and $c$ does not depend on $\omega$ as long as $\omega\leq\omega_+$.
The trace theorem allows us to estimate the coefficients $a_m$ in terms of $u$ in $\Omega$,
\begin{equation}\label{esta}
  \sum_{m\in\zprop} |a_m^\text{inc}+a_m|^2 (1+|m|) + \sum_{m\not\in\zprop} |a_m|^2 (1+|m|) = \|u|_{\Gamma_-}\|^2_{H^{1/2}(\Gamma)} \leq M\Hnorm{1}{\Omega}{u}^2.
\end{equation}
From estimates (\ref{estu1},\ref{estu},\ref{esta}), we obtain
\begin{equation}
  c\Hnorm{1}{\tilde\Omega}{u}^2 \leq M\Hnorm{1}{\Omega}{u}^2 + 2\!\!\sum_{m\in\zprop} |a_m^\text{inc}|^2 (1+|m|).
\end{equation}
A similar estimate is obtained for a domain analogous to $\tilde\Omega$ for $x_3>z_+$.  As a result of the pseudoperiodicity of $u$ and the boundedness of $D$, we obtain the desired estimate.
\end{proof}

\subsection{Eigenvalues of the scattering problem} 

We will assume that the functions $\eps(x)$ and $\tau(x)$ are bounded from below and above in $\Omega$ by fixed positive constants,
\begin{equation}\label{bounds}
        0<\eps_-^0 \leq \eps(x) \leq \eps_+^0
       \quad\text{and}\quad
       0<\tau_-^0 \leq \tau(x) \leq \tau_+^0
       \quad\text{for $x\in\Omega$.}
\end{equation}

The weak form of the scattering problem can be expressed in terms of the following sesquilinear forms in $\Honek$:
\begin{eqnarray}
  && a^\omega(u,v) = \int_\Omega
     \tau\, \nabla u \cdot \nabla \bar v\; dV
     + \tau_0 \int_\Gamma (T^\omega u) \bar v \; dA, \\
  && a^\omega_r(u,v) = \int_\Omega
     \tau\, \nabla u \cdot \nabla \bar v\; dV
     + \tau_0 \int_\Gamma (T^\omega_r u) \bar v \; dA, \\
   && a^\omega_i(u,v) = \tau_0 \int_\Gamma (T^\omega_i u) \bar v \; dA, \\
   && b(u,v) = \int_\Omega \eps\, u \bar v \; dV.
\end{eqnarray}
Observe that $a^\omega = a^\omega_r + ia^\omega_i$.
In terms of these forms, Problem \ref{problemscattering2} can be written as
 \begin{equation}\label{inhomogeneous}
  a^\omega(u,v) - \omega^2 b(u,v) = f^\omega_{\Gamma}(v)  \quad \text{for all $v\in\Honek$,}
\end{equation}
and the generalized scattering problem as
\begin{equation}\label{generalized}
  a^\omega(u,v) - \omega^2 b(u,v) = f(v)  \quad \text{for all $v\in\Honek$} \quad (f\in\Honek^*).
\end{equation}
We consider first the homogeneous problem
\begin{equation}\label{homogeneous}
  a^\omega(u,v) - \omega^2 b(u,v) = 0 \quad \text{for all $v\in\Honek$.}
\end{equation}
This is a nonlinear eigenvalue problem because of the dependence of $a^\omega$ on $\omega$ through the Dirichlet-to-Neumann operator~$T^\omega$.

\begin{definition}
 A number $\omega$ is said to be an eigenvalue of a one-parameter family of bounded sesquilinear forms $c^\omega(\cdot,\cdot)$ in a Hilbert space $H$ if there exists a nonzero element $u\in H$ such that, for all $v\in H$, $c^\omega(u,v)=0$.
\end{definition}

The eigenvalues of the family $a^\omega - \omega^2b$ are in general complex.  Its real eigenvalues form a subset of the eigenvalues of  the real part of the form, namely $a_r^\omega - \omega^2b$, as stated in Proposition~\ref{propcharacterization} below.  For an eigenfunction of the real form to be an eigenvalue of the complex form also, all of its propagating Fourier harmonics must vanish.  This means that, as long as $\zlin$ is empty, a nontrivial solution of \eqref{homogeneous}, for real $\omega^2$, falls off exponentially with distance from the slab structure; such a field is a {\em guided mode} of the slab.  If this frequency $\omega$ is large enough so that $\zprop$ is not empty, then $\omega$ is an embedded eigenvalue for the $\kappa$-pseudoperiodic operator corresponding to the partial (in $x'=(x_1,x_2)$) Floquet-Bloch decomposition of the Helmholtz equation in~$\RR^3$.  Typically, an embedded eigenvalue is not robust with respect to perturbations of $\kappa$, $\eps$, or $\tau$ because the condition \eqref{char2} below that the coefficients of all propagating harmonics vanish is generically not satisfied.  The existence of a guided mode requires special conditions, such as symmetry of $\eps$ and $\tau$ for $\kappa=0$.  The reader is referred to \cite{Bonnet-BeStarling1994}, \cite{ShipmanVolkov2007}, and \cite{ShipmanVenakides2005} for further discussion of non-robust guided modes in this context.

\begin{proposition}[Characterization of real eigenvalues]\label{propcharacterization} \hspace{2ex}
  If $\omega^2~\in~\RR$, then a function $u\in\Honek$ satisfies the homogeneous problem \eqref{homogeneous} if and only if it satisfies the equation
\begin{equation}\label{adjointform}
  a_r^\omega(u,v) - ia_i^\omega(u,v) - \omega^2b(u,v) = 0 \quad \text{for all $v\in\Honek$}
\end{equation}
and if and only if it satisfies the pair
\begin{eqnarray}
  && a^\omega_r(u,v) - \omega^2 b(u,v) = 0 \quad \text{for all $v\in\Honek$,} \label{char1} \\
  && (\widehat{u|_\Gamma})_m = 0 \quad \text{for all $m\in\zprop$}. \label{char2}
\end{eqnarray}
\end{proposition}

\begin{proof}
We prove that \eqref{homogeneous} is equivalent to the pair (\ref{char1},\ref{char2}).  The equivalence to equation \eqref{adjointform} is proved similarly.  Suppose that $\omega$ and $u\not=0$ satisfy \eqref{homogeneous}.  The imaginary part of this equation with $v=u$, together with the expression \eqref{Ti} for $T_i$, gives
\begin{equation}
  \tau_0 \sum_{m\in\zprop} -\eta_m|(\widehat{u|_\Gamma})_m|^2 = \int_\Gamma (T_i^\omega u)\bar u = 0.
\end{equation}
Since $\eta_m>0$ for all $m\in\zprop$, all propagating Fourier coefficients $(\widehat{u|_\Gamma})_m$ of $u$ on $\Gamma$ vanish.  This in turn proves that $a_i(u,v)=0$ for all $v\in\Honek$, so that $\omega$ and $u$ satisfy \eqref{char1}.
Conversely, if \eqref{char2} holds, then $a_i(u,v)=0$ for all $v\in\Honek$ and therefore \eqref{char1} is equivalent to \eqref{homogeneous}.
\end{proof}

\begin{proposition}[Real eigenvalue sequences]\label{propeigenvalues}
  Given the bounds \eqref{bounds} on the functions $\eps$ and $\tau$, the eigenvalues $\omega$ of the family $a^\omega_r(u,v) - \omega^2 b(u,v)$ consist of the elements of a nondecreasing sequence of positive numbers $\{ \omega_j(\eps,\tau) \}_{j=1}^\infty$ that tends to $\infty$ and their additive inverses.  The eigenvalues of the family $a^\omega(u,v)-\omega^2b(u,v)$ consist of a subsequence $\{ \omega^*_k(\eps,\tau) \}_{k=1}^N$ of this sequence, where $N$ is a nonnegative integer (perhaps 0) or infinity.
\end{proposition}

\begin{proof}
The proof follows \cite{Bonnet-BeStarling1994}.  As the family $a^\omega_r(u,v) - \omega^2 b(u,v)$ depends only on $\omega^2$, we shall consider only nonnegative values of $\omega$.  Let $\omega\geq0$ be given, and let us consider the set of numbers $\lambda$ such that there exists a nonzero function $u\in\Honek$ such that $a^\omega_r(u,\cdot) - \lambda b(u,\cdot) = 0$.  According to the {\em min-max principle} (see \cite{ReedSimon1980d}, \S XIII, for example), this set consists of a strictly increasing sequence of positive numbers $\{\lambda^\omega_j(\eps,\tau)\}_{j=1}^\infty$ defined by
\begin{multline}\label{Rayleigh}
  \lambda^\omega_j(\eps,\tau) = \sup_{V^{j-1}<L^2(\Omega)}\inf_{\mbox{\parbox{14ex}{\scriptsize $u\in(V^{j-1})^\perp\!\setminus\!\{0\}$ \\ \hspace*{4ex}$u\in\Honek$}}}
  \frac{a^\omega_r(u,u)}{b(u,u)} = \\
 = \sup_{V^{j-1}<L^2(\Omega)}\inf_{\mbox{\parbox{14ex}{\scriptsize $u\in(V^{j-1})^\perp\!\setminus\!\{0\}$ \\ \hspace*{4ex}$u\in\Honek$}}} \frac{\int_\Omega\tau\left| \nabla u \right|^2\; dV + \tau_0\int_\Gamma\left( T^\omega_r u \right)\bar u\; dA}{\int_\Omega \eps |u|^2 \; dA},
\end{multline}
in which the supremum is taken over all $k$-dimensional subspaces $V^k$ of $L^2(\Omega)$, for $k~=~j~-~1$, and ``$\perp$" refers to the orthogonal complement with respect to the norm $b(u,u)$ in $L^2(\Omega)$.  One can prove that, for each positive integer $j$, $\lambda_j^\omega(\eps,\tau)$ is a continuous and nonincreasing function of $\omega\!\geq\!0$ (see the proof of Theorem 3.3 of \cite{Bonnet-BeStarling1994}).  There is therefore, for each $j$, exactly one positive number, which we denote by $\omega_j(\eps,\tau)$, that satisfies
\begin{equation}
\omega_j(\eps,\tau)^2 = \lambda_j^{\omega_j(\eps,\tau)}(\eps,\tau).
\end{equation}
The number $\omega$ is an eigenvalue of the family $a^\omega_r(u,v) - \omega^2 b(u,v)$ if and only if there exists an integer $j$ such that $\omega^2 = \lambda^\omega_j(\eps,\tau)$.  The sequence $\{ \omega_j(\eps,\tau) \}_{j=1}^\infty$ therefore consists of all the nonnegative eigenvalues of the family.

The second statement in the Proposition follows from this and Proposition \ref{propcharacterization}.
\end{proof}


\medskip

Because the Rayleigh quotient in the min-max principle \eqref{Rayleigh} decreases with an increase in $\eps$ and increases with an increase in $\tau$, the eigenvalues inherit the property of monotonicity with respect to these functions.

\begin{proposition}[Eigenvalue dependence on $\eps$ and $\tau$]\label{propmonotonicity}
Let $\eps_-$, $\eps_+$, $\tau_-$, and $\tau_+$ be measurable real-valued functions on $\Omega$ that satisfy the bounds \eqref{bounds} and the inequalities $\eps_-(x)\leq\eps_+(x)$ and $\tau_-(x)\leq\tau_+(x)$ on $\Omega$.  Then, for each positive integer~$j$,
\begin{equation}
  \omega_j(\eps_+,\tau_-) \leq \omega_j(\eps_-,\tau_+).
\end{equation}
\end{proposition}

\section{Proof of the main theorem}\label{secmain} 

Proof of the theorem on differentiability of the solution of $u$ and the transmitted energy~$\trans$ with respect to $\eps$ and $\tau$ rely on Meyers' theorem on higher regularity of solutions of elliptic equations.  As we have discussed,
in order to apply this theorem to the solution $u$ of the scattering problem, it is necessary to be assured that $u$ is uniformly bounded over all admissible functions $\eps$ and $\tau$.
The precise condition we will need is one on lower and upper bounding functions for these material coefficients.
 
\begin{condition}[Non-resonance]\label{condnonresonance}
For a given number $\omega\in\RR$, the measurable real-valued functions $\eps_-$, $\eps_+$, $\tau_-$, and $\tau_+$ on $\Omega$ satisfy the {\em non-resonance condition} if, for each pair ($\eps,\tau$) of measurable real-valued functions on $\Omega$ that satisfy
\begin{equation}
  \eps_-(x) \leq \eps(x) \leq \eps_+(x)
  \quad\text{and}\quad
  \tau_-(x) \leq \tau(x) \leq \tau_+(x),
\end{equation}
for all $x\in\Omega$, $\omega$ is not an eigenvalue of the family $a^\omega - \omega^2b$.
\end{condition}

This condition can be arranged if we choose the upper and lower bounding functions such that, for some integer $j$,
\begin{equation}\label{evalcondition}
  \omega_j(\eps_-,\tau_+) < \omega_{j+1}(\eps_+,\tau_-).
\end{equation}
Then, Proposition \ref{propmonotonicity} guarantees that, for all functions $\eps$ and $\tau$ between these functions,
\begin{equation}
  \omega_j(\eps,\tau) \leq \omega_j(\eps_-,\tau_+) < \omega_{j+1}(\eps_+,\tau_-) \leq \omega_{j+1}(\eps,\tau),
\end{equation}
in which case Condition \ref{condnonresonance} holds for each $\omega$ strictly between $\omega_j(\eps_-,\tau_+)$ and $\omega_{j+1}(\eps_+,\tau_-)$.  The condition \eqref{evalcondition} can be achieved, for example, by beginning with a fixed pair of functions $(\eps,\tau)$ for which $\omega_j(\eps,\tau) < \omega_{j+1}(\eps,\tau)$, and varying them up and down continuously in the $L^\infty$ norm, with respect to which each $\omega_j(\eps,\tau)$ is a continuous function of $\eps$ and $\tau$.
In view of the fact that the eigenvalues of $a^\omega -\omega^2b$ are a subset of $\{\omega_j\}_{j=1}^\infty$, the condition \eqref{evalcondition} is evidently stronger than what is necessary.  In fact, as we have discussed, any $\omega_j$ is typically not an eigenvalue of the scattering problem, for given material coefficients $\eps$ and $\tau$.

\subsection{A uniform bound for the scattered field}

The following theorem guarantees a bound on the solution of the scattering problem that is uniform over functions $\eps$ and $\tau$ bounded below and above by functions that satisfy Condition~\ref{condnonresonance}.

\begin{theorem}[Bound on the scattered field]\label{thmfieldbound}
Let $\eps_-$, $\eps_+$, $\tau_-$, and $\tau_+$ be measurable real-valued functions on $\Omega$ that satisfy the bounds \eqref{bounds} and the non-resonance Condition \ref{condnonresonance} for all $\omega$ in some positive interval $[\omega_-,\omega_+]$.  There exists a positive number $K$ such that, for each $\omega\in[\omega_-,\omega_+]$ and each pair of measurable real-valued functions $\eps$ and $\tau$ on $\Omega$ that satisfy
\begin{equation}\label{bounds2}
  \eps_-(x) \leq \eps(x) \leq \eps_+(x)
  \quad\text{and}\quad
  \tau_-(x) \leq \tau(x) \leq \tau_+(x),
\end{equation}
the generalized scattering problem \eqref{weakgeneral} admits a unique solution $u$ such that
\begin{equation}
  \Hnorm{1}{\Omega}{u}<K\| f \|_{\Honek^*}, 
\end{equation}
where $f\in\Honek^*$ denotes the general functional on the right-hand side of \eqref{weakgeneral}.
\end{theorem}

\begin{proof}
We first prove that the scattering problem \eqref{weakgeneral} admits a unique solution for the parameters given in the Theorem.  Rewrite \eqref{generalized}~as
\begin{equation}\label{rewrit}
  [a^\omega(u,v) + b(u,v)] - (\omega^2+1)b(u,v)  =  f(v).
\end{equation}
Since both $a$ and $b$ are bounded forms in $\Honek$, there exist linear operators $A^\omega$ and $C^\omega$ from $\Honek$ into itself, as well as an element $\tilde f \in\Honek$ defined through
\begin{eqnarray}
  && (A^\omega u,v) = a^\omega(u,v) + b(u,v), \\
  && (C^\omega u,v) = - (\omega^2+1)b(u,v), \\
  && (\tilde f,v) = f(v).
\end{eqnarray}
In terms of these objects, equation \eqref{rewrit} takes the form
\begin{equation}\label{inhomog}
  (A^\omega + C^\omega) u = \tilde f.
\end{equation}
The operator $A^\omega$ is bijective with a bounded inverse because $a^\omega(u,v) + b(u,v)$ is coercive (recall that $T^\omega_r$ is a positive operator):
\begin{multline}
  \Re ( a^\omega(u,u) + b(u,u) ) = \int_\Omega \tau |\nabla u |^2 + \int_\Omega \eps |u|^2
                + \tau_0 \int_\Gamma (T^\omega_r u)\bar u \\
  \geq \tau_-^0 \int_\Omega |\nabla u|^2 + \eps_-^0 \int_\Omega |u|^2
  \geq \min \left\{ \tau_-^0, \eps_-^0 \right\} \| u \|_{\Honek}.
\end{multline}
Moreover, $C^\omega$ is compact because of the compact embedding of $\Honek$ into $L^2(\Omega)$.
By the Fredholm alternative, \eqref{inhomog} (equivalently, \eqref{rewrit}) has a unique solution if $A^\omega+C^\omega$ is injective, that is, $\omega$ is not an eigenvalue of the family $a^\omega-\omega^2b$.  But this is implied by the non-resonance Condition \ref{condnonresonance} which we have assumed for $\omega$.

We turn to establishing a bound on this solution that is uniform over all functions $\eps$ and $\tau$ and numbers $\omega$ that satisfy the hypotheses of the Theorem and all $f$ in the unit ball in $\Honek^*$.
To accomplish this, it suffices to consider arbitrary sequences $\eps_n$ and $\tau_n$ of measurable functions that satisfy the bounds \eqref{bounds2}, a sequence $\omega_n$ of numbers satisfying $\omega_-\leq \omega_n \leq \omega_+$, and sequences $u_n\in\Honek$ and $f_n\in\Honek^*$ with $\Hnorm{1}{\Omega}{u_n}\leq1$ and $f_n\to0$ such that
\begin{equation}
  \int_\Omega
\left( \tau_n\, \nabla u_n \!\cdot\! \nabla \bar v - \omega_n^2\,\eps_n  u_n \bar v  \right)
+ \tau_0 \int_\Gamma (T^{\omega_n} u_n) \bar v
\,=\, f_n(v)
\quad
\text{for all $v\in\Honek$},
\end{equation}
and to prove that, necessarily, $\Hnorm{1}{\Omega}{u_n}\to0$.
We may as well assume (by extracting a subsequence) that there exists a number $\omega\in[\omega_-,\omega_+]$ such that $\omega_n\to\omega$.  We rewrite this equation as
\begin{equation}\label{eqng}
  \int_\Omega
\left( \tau_n\, \nabla u_n \!\cdot\! \nabla \bar v - \omega^2\,\eps_n u_n \bar v  \right)
+ \tau_0 \int_\Gamma (T^{\omega} u_n) \bar v
\,=\, g_n(v)
\quad
\text{for all $v\in\Honek$},
\end{equation}
in which the elements $g_n\in\Honek^*$ are defined by
\begin{equation}\label{defng}
  g_n(v) = f_n(v) + (\omega^2-\omega_n^2)\!\int_\Omega \eps_n u_n \bar v
       + \tau_0 \int_\Gamma (T^{\omega} - T^{\omega_n}) u_n \bar v.
\end{equation}
We shall prove that $g_n\to0$ in $\Honek^*$.

We first estimate the third term in \eqref{defng},
\begin{equation}
  \tau_0 \int_\Gamma(T^\omega-T^{\omega_n})u_n\bar v = -i\tau_0 \sum_{m\in\ZZ^2} (\eta^\omega_m-\eta^{\omega_n}_m)\hat u_{nm}\bar{\hat v}_m,
\end{equation}
in which, for simplicity, $\hat u_{nm}$ denotes the $\kappa$-Fourier coefficient of $u_{mn}|_\Gamma$.
It is straightforward from the definition of $\eta_m$ to demonstrate that there exists a number $c$ such that, for all $m\in\ZZ^2$ and all $n$ sufficiently large,
\begin{equation}
  | \eta^\omega_m - \eta^{\omega_n}_m | \leq c | \omega^2 - \omega_n^2 |^{1/2}.
\end{equation}
This allows us to estimate, using $\| u_n \|_{\Honek}\leq 1$,
\begin{equation}
  \left| -i\tau_0 \sum_{m\in\ZZ^2} (\eta^\omega_m-\eta^{\omega_n}_m)\hat u_{nm}\bar{\hat v}_m \right|
  \leq \tau_0 c | \omega^2 - \omega_n^2|^{1/2} \| u_n |_\Gamma \|_{L^2(\Omega)} \| v |_\Gamma \|_{L^2(\Omega)}
  \leq \tau_0 M^2 c | \omega^2 - \omega_n^2 |^{1/2} \| v \|_{\Honek}.
\end{equation}
This proves that
\begin{equation}\label{gconv1}
  \tau_0  \int_\Gamma (T^\omega - T^{\omega_n}) u_n \bar v \to 0 \text{ as } n\to\infty
\end{equation}
uniformly over $v\in\Honek$ with $\| v \| = 1$.

The second term of \eqref{defng} is estimated~by
\begin{equation}
  \left| (\omega^2 - \omega_n^2) \int_\Omega \eps_n u_n \bar v \right|
  \leq  \eps^0_+ | \omega^2 - \omega_n^2 | \| u_n \|_{L^2(\Omega)} \| v \|_{L^2(\Omega)}
  \leq M^2 \eps^0_+ | \omega^2 - \omega_n^2 | \| v \|_{\Honek},
\end{equation}
demonstrating that
\begin{equation}\label{gconv2}
  (\omega^2 - \omega_n^2) \int_\Omega \eps_n u_n \bar v \to 0 \text{ as } n\to \infty
\end{equation}
uniformly over $v\in\Honek$ with $\| v \| = 1$.

The results (\ref{gconv1},\ref{gconv2}), together with the convergence $f_n\to0$ prove the strong convergence of $g_n$ to zero in $\Honek^*$.

We shall now demonstrate that there exists a function $u\in\Honek$, measurable functions $\eps$ and $\tau$ satisfying the bounds \eqref{bounds2}, and an infinite subset $\Upsilon$ of the positive integers such that the following convergences hold, restricted to indices in the subsequence $\Upsilon$,
\begin{equation}\label{convergencelist}
\renewcommand{\arraystretch}{1.5}
\left.
  \begin{array}{ll}
      u_n \rightharpoonup u & \text{weak } H^1(\Omega), \\
      u_n \to u & \text{strong } L^2(\Omega), \\
      \eps_n \rightharpoonup \eps & \text{weak* } L^\infty(\Omega), \\
      \eps_n u_n \to \eps u & \text{strong } H^{-1}(\Omega), \\
      \displaystyle \tau_n \to {\tau} & \text{G-convergence in } \Omega, \\
      g_n \to 0 & \text{strong } H^{-1}(\Omega).
  \end{array}
\right.
\end{equation}
The first and second subsequence limits are due to the uniform bound on the functions $u_n$ in $\Hone$, the Alaoglu Theorem, and the compact embedding of $\Hone$ into $\Ltwo$.  The third is due to the uniform bound on the functions $\eps_n$ in $L^\infty(\Omega)$ and the Alaoglu Theorem.  Because of the strong $L^2$ convergence of $u_n$, we obtain, for each $v\in\Ltwo$, $u_n\bar v \to u\bar v$ in $L^1(\Omega)$ (for the subsequence $\Upsilon$), and therefore because of the weak-* convergence of $\eps_n$,
\begin{equation}
  \int_\Omega \eps_n u_n \bar v \to \int_\Omega \eps u\bar v
  \quad
  \text{for all } v\in\Honek \supset \Honez,
\end{equation}
from which we infer that $\eps_n u_n\rightharpoonup\eps u$ weakly in $H^{-1}(\Omega)=\Honez^*$, or, more precisely, that $j(\eps_n u_n)\rightharpoonup j(\eps u)$, where $j$ is the natural embedding of $\Ltwo$ into $H^{-1}(\Omega)$ defined by $j(w)(v) = \int_\Omega w\bar v$ for $v\!\in\!\Honez$.  Since this embedding is compact and the sequence $\eps_n u_n$ is bounded in $\Ltwo$, we have strong (and therefore also weak) convergence of a subsequence, $\eps_n u_n \to w\in H^{-1}(\Omega)$, and by the uniqueness of weak limits, we obtain $w=\eps u$, which justifies the fourth convergence in the list \eqref{convergencelist}.  The last convergence follows from the strong convergence $g_n\to0$ in $\Honek^*$ and the inclusion $\Honez\subset\Honek$.
The existence of the $G$-limit (or $H$-limit) $\tau$ satisfying the bounds $\tau_-(x)\leq\tau(x)\leq\tau_+(x)$ for $x\in\Omega$ follows from Theorem 2 of Murat/Tartar \cite{MuratTartar1978} and the discussion in the second paragraph of that work (p. 21).  

The divergence $\Div\xi$ of a vector field $\xi\in\Ltwo$ is the element of $H^{-1}(\Omega)$ defined by
\begin{equation}
  (\Div\xi)(v) = \int_{\Omega} \xi\cdot\nabla v
  \quad
  \text{for } v\in\Honez,
\end{equation}
whose norm is bounded by the norm of $\xi$,
\begin{equation}
  \| \Div \xi \|_{H^{-1}(\Omega)} \leq \| \xi \|_{\Ltwo}.
\end{equation}
From equation \eqref{eqng} and items 4 and 6 in \eqref{convergencelist}, we infer
\begin{equation}\label{divconvergence}
  \Div \tau_n \nabla u_n \to j(\eps u) \quad \text{strongly in } H^{-1}(\Omega)
\end{equation}
(the action of the integral over $\Gamma$ in \eqref{eqng} is trivial on $\Honez$).  Because of the strong convergence of $\Div \tau_n \nabla u_n$, the weak convergence of $u_n$ in $\Hone$ and the G-convergence of $\tau_n$, we may apply Theorem 1 of \cite{MuratTartar1978} to deduce that
\begin{equation}
  \tau_n \nabla u_n \rightharpoonup \tau \nabla u \quad \text{weakly in } \Ltwo.
\end{equation}
Because of the weak convergence of $u_n$ in $\Hone$, we have, for all $v\in\Honek$ (for $n\in\Upsilon$)
\begin{equation}
  \tau_0 \int_\Gamma (T^\omega u_n) \bar v = \tau_0 \int_\Gamma u_n (\overline{{T^\omega}^* v})
  \longrightarrow
  \tau_0 \int_\Gamma u (\overline{{T^\omega}^* v}) = \tau_0 \int_\Gamma (T^\omega u) \bar v.
\end{equation}
We can now take the limit of each term in \eqref{eqng} to obtain
\begin{equation}\label{limiteqn}
  \int_\Omega
\left( \tau\, \nabla u \!\cdot\! \nabla \bar v - \omega^2\,\eps \, u \bar v  \right)
+ \tau_0 \int_\Gamma (T^{\omega} u) \bar v
\,=\, 0
\quad
\text{for all $v\in\Honek$}.
\end{equation}
By the uniqueness of the solution to this problem, which we proved above, we must have $u=0$ in $\Hone$.
Equation \eqref{eqng}, with $v$ set equal to $u_n$, gives
\begin{multline}
  \eps_+^0\omega^2\Lnorm{2}{\Omega}{u_n}^2 \geq \omega^2 \int_\Omega \eps_n|u_n|^2
  \geq \tau_-^0 \int_\Omega |\nabla u_n|^2 - |g_n(u_n)| \\
  = \tau_-^0 (\Hnorm{1}{\Omega}{u_n}^2 - \Lnorm{2}{\Omega}{u_n}^2) - |g_n(u_n)|,
\end{multline}
whence we obtain
\begin{equation}
  (\eps_+^0\omega^2 + \tau_-^0)\Lnorm{2}{\Omega}{u_n}^2
  \geq \tau_-^0 \Hnorm{1}{\Omega}{u_n}^2 - |g_n(u_n)|.
\end{equation}
From of the strong convergence $u_n\!\to\!u=0$ in $\Ltwo$ and the strong convergence of $g_n\to0$ in $\Honek^*$, we deduce that $\Hnorm{1}{\Omega}{u_n}\to0$, as we set out to do.  We conclude that there exists a number $K$ such that the solution $u$ of the generalized scattering problem \eqref{weakgeneral} satisfies
\begin{equation}
  \| u \|_{\Honek} \leq K \| f \|_{\Honek^*}
\end{equation}
for all $f\in\Honek^*$, for all functions $\eps$ and $\tau$ that satisfy \eqref{bounds2}, and for all $\omega\in[\omega_-,\omega_+]$.
\end{proof}
\bigskip

\subsection{Field sensitivity to $L^p$ perturbations}

This section contains the main theorem of this work, Theorem~\ref{thmmain}, and its proof.  The theorem makes rigorous the formal variational gradient, obtained in section \ref{secformal}, of the solution $u$ of the scattering problem as a function of the material coefficients $\eps$ and $\tau$.  The field $u$ satisfies
\begin{equation}\label{weak1}
\int_\Omega
\left( \tau \, \nabla u \cdot \nabla \bar v - \omega^2\, \eps\, u \bar v  \right)
+ \tau_0 \int_{\Gamma} (T^\omega u) \bar v
\,=\, \tau_0 \int_{\Gamma} \left( \partial_n + T^\omega \right)\inc \bar v
\quad\text{for all $v\in\Honek$}.
\end{equation}
If we replace $\eps$, $\tau$, and $u$ with $\eps+\breve\eps$, $\tau+\breve\tau$, and $u+\breve u$ and subtract from \eqref{weak1}, we obtain
\begin{equation}\label{perturbed}
\int_\Omega
\left( \tau \, \nabla \breve u \cdot \nabla \bar v - \omega^2\, \eps\, \breve u \bar v  \right)
+ \tau_0 \int_{\Gamma} (T^\omega \breve u) \bar v
\,=\, 
- \int_\Omega
\left( \breve \tau \, \nabla (u+\breve u) \cdot \nabla \bar v - \omega^2\, \breve\eps\, (u+\breve u) \bar v  \right).
\end{equation}
Retaining only the linear part of the right-hand side gives an equation for $\breve u_0$, the formal linearization of the perturbation of $u$ about $(\eps,\tau)$,
\begin{equation}\label{approx}
\int_\Omega
\left( \tau \, \nabla \breve u_0 \cdot \nabla \bar v - \omega^2\, \eps\, \breve u_0 \bar v  \right)
+ \tau_0 \int_{\Gamma} (T^\omega \breve u_0) \bar v
\,=\, 
- \int_\Omega
\left( \breve \tau \, \nabla u \cdot \nabla \bar v - \omega^2\, \breve\eps\, u \bar v  \right).
\end{equation}
The task is to prove that $\|\breve u - \breve u_0\| = \littleo(\|(\breve\eps,\breve\tau)\|)$ as $\breve\eps$ and $\breve\tau$ tend to zero in an $L^p$ norm.

Recall that, if, for some vector function $\xi\in\Ltwo$ and scalar function $h\in\Ltwo$, a function $u\in\Honez$ satisfies
\begin{equation}
\int_\Omega
\left( \tau \, \nabla u \cdot \nabla \bar v - \omega^2\, \eps\, u \bar v  \right)
\,=\, 
\int_\Omega
\left( \xi \cdot \nabla \bar v + h \bar v  \right)
\end{equation}
for all $v\in\Honez$, we say that $u$ satisfies, in the weak sense, the partial differential equation
\begin{equation}\label{pde}
  \Div \tau\nabla u - \omega^2 \eps u = \Div \xi + h  \quad \text{in } \Omega.
\end{equation}

The proof of Theorem~\ref{thmmain} requires the following specialization of the theorem of Meyers
on the higher integral regularity of solutions of elliptic differential equations.\footnote{
In Theorem 2 of \cite{Meyers1975}, we fix $p_1=2$, $r=2$, and the dimension $n=3$.  The $p$ in Meyers' theorem corresponds to our $q$ here.  We also enforce $Q\leq6$, which guarantees $r^*\geq q$ for all $q<Q$, because
$r^* = (r^{-1}-n^{-1})^{-1} = 6 \geq Q > q$.  We may use equation (49) from the theorem because $q>2>6/5 = 2n/(n+2)$ with $n=3$.}

\begin{theorem}[Meyers regularity]\label{thmMeyers}

Given a bounded domain $D\subset\RR^3$, a real-valued measurable function $\tau(x)$ in $D$, and positive real numbers $\tau_-$, $\tau_+$, and $R_0$ such that
\begin{equation}
  0 < \tau_- \leq \tau(x) \leq \tau_+,
\end{equation}
there exists a number $Q$ with $2<Q<6$ such that, for each $q$ satisfying
\begin{equation}
  2<q<Q<6,
\end{equation}
there exist constants $C_1$ and $C_2$ such that, given
\begin{eqnarray}
  && u\in H^1(D), \; \xi\in L^q(D), \; h\in L^2(D), \\
  && \nabla \cdot \tau \nabla u = \Div \xi + h \text{ in $D$}, \\
  && B_y(2R)\subset D, \; \text{with $y\in D$ and $R>R_0$},
\end{eqnarray}
the following inequalities hold:
\begin{eqnarray}
  \Lnorm{q}{B_y(R)}{\nabla u} &<&  \label{Meyers1}
    C_1\left[ \Lnorm{2}{B_y(2R)}{u} + \Lnorm{q}{B_y(2R)}{\xi} + \Lnorm{2}{B_y(2R)}{h}\right], \\
  \Lnorm{q}{B_y(R)}{u} &<&
    C_2\Hnorm{1}{B_y(R)}{u}.   \label{Meyers2}
\end{eqnarray}
\end{theorem}

The second statement \eqref{Meyers2} is a result of the compact embedding of $H^1(B_y(R))$ into $L^q(B_y(R))$ for $q<6$ (see Theorem 7.26 of \cite{GilbargTrudinger1998}, for example).

\medskip

\begin{theorem}[Field sensitivity to $L^p$ variations]\label{thmmain}
Let $\eps_-$, $\eps_+$, $\tau_-$, and $\tau_+$ be measurable real-valued functions on $\Omega$ that satisfy the bounds \eqref{bounds} and the non-resonance Condition \ref{condnonresonance} for all $\omega$ in some positive interval $[\omega_-,\omega_+]$.  Assume additionally that $\sum_{m\in\zprop}(|a_m^\text{\rm inc}|^2 + |b_m^\text{\rm inc}|^2)$ is bounded uniformly for $\omega\in[\omega_-,\omega_+]$.
Then there exist real numbers $C$ and $p\!>\!6$ such that, for all $\omega\in[\omega_-,\omega_+]$ and all
measurable functions $\eps$, $\breve\eps$, $\tau$, and $\breve\tau$ on $\Omega$ that satisfy
\begin{eqnarray}\label{bounds3}
  \eps_-(x) \leq \eps(x) \leq \eps_+(x)
     &\text{and}&
     \tau_-(x) \leq \tau(x) \leq \tau_+(x),  \\ \label{bounds4}
  \eps_-(x) \leq (\eps+\breve\eps)(x) \leq \eps_+(x)
     &\text{and}&
     \tau_-(x) \leq (\tau+\breve\tau)(x) \leq \tau_+(x),
\end{eqnarray}
the following statement holds:

If $u\in\Honek$ is the unique solution of the scattering problem guaranteed by Theorem~\ref{thmfieldbound} (that is, $u$ satisfies \eqref{weak1} for all $v\in\Honek$), $u+\breve u$ is the unique solution of the scattering problem with $\tau$ replaced by $\tau+\breve\tau$ and $\eps$ replaced by $\eps+\breve\eps$ in \eqref{weak1}, and
$\uz$ satisfies the approximate equation \eqref{approx} for all $v\in\Honek$, then the linear operator
\begin{equation}\label{derivative}
  \left( L^p(\Omega) \right)^2\to H^1(\Omega) :: (\breve\tau,\breve\eps)\mapsto \breve u_0 
\end{equation}
(restricted to $(\beps,\btau)$ admissible by \eqref{bounds4}) satisfies
\begin{equation}
  \Hnorm{1}{\Omega}{\breve u_0} \leq C \left( \Lnorm{p}{\Omega}{\breve\tau} + \Lnorm{p}{\Omega}{\breve\eps} \right)
\end{equation}
and
\begin{equation}
  \Hnorm{1}{\Omega}{\breve u - \uz} \leq C \left( \Lnorm{p}{\Omega}{\breve\tau} + \Lnorm{p}{\Omega}{\breve\eps} \right)^2.
\end{equation}

Moreover, the derivative \eqref{derivative} is Lipschitz continuous, that is, if
$(\btau,\beps)\mapsto \buz^1$ and $(\btau,\beps)\mapsto\buz^2$ denote the derivatives at $(\eps^1,\tau^1)$ and $(\eps^2,\tau^2)$, respectively, then
\begin{equation}\label{Lipschitz}
  \Hnorm{1}{\Omega}{\buz^1-\buz^2} < C \left( \Lnorm{p}{\Omega}{\tau^1\!-\!\tau^2} + \Lnorm{p}{\Omega}{\eps^1\!-\!\eps^2} \right) \left( \Lnorm{p}{\Omega}{\btau} + \Lnorm{p}{\Omega}{\beps} \right),
\end{equation}
so long as the functions $\eps^1$, $\eps^2$, $\tau^1$, $\tau^2$, $\eps^1+\beps$, $\eps^2+\beps$, $\tau^1+\btau$, and $\tau^2+\btau$ satisfy the bounds (\ref{bounds3},\ref{bounds4}).
\end{theorem}

Because the conclusion of the theorem holds for each $q\geq p$ if it holds for a given $p$, the condition $p>6$ could be logically be replaced, equivalently, with the condition $p\geq1$.  We have used the number 6 because, in the proof, $p$ arises as a number greater than 6.

\bigskip
\begin{proof} 
Let $\eps$, $\breve\eps$, $\tau$, and $\breve\tau$ satisfy the bounds in the Theorem, and let $\omega\in[\omega_-,\omega_+]$ be given. 
Let $B_0$ be a ball of radius $R_0$ containing $\Omega$, and let $B_1$ and $B_2$ be the balls whose centers coincide with that of $B_0$ and whose radii are $2R_0$ and $4R_0$, respectively.
Let $Q$ be as provided in Theorem~\ref{thmMeyers} for the domain $B_2$ and the constant $R_0$.
Let $s$, $q$, and $p$ be such that
\begin{equation}
  2 < s < q < Q,
\end{equation}
\begin{equation}
  q^{-1} + p^{-1} = s^{-1} \quad \text{and} \quad s^{-1} + p^{-1} = 2^{-1}.
\end{equation}
Let the constants $C_1$ and $C_2$ in Theorem~\ref{thmMeyers} be valid for both $q$ and $s$ in place of the $q$
in the Theorem.  Denote the solutions of the scattering problems in $\RR^3$ corresponding to $u$ and $u+\breve u$ by these same symbols.

Theorem~\ref{thmfieldbound} and Lemma \ref{lemmaextension} together provide a number $K_1$, independent of the choice of $\eps$, $\breve\eps$, $\tau$, and $\breve\tau$, such that
\begin{equation}\label{uniformbound}
  \Hnorm{1}{B_2}{u+\breve u} < K_1 \quad\text{and}\quad
  \Hnorm{1}{B_2}{u} < K_1.
\end{equation}
Because the condition on the incident field stated in the Theorem makes $f^\omega_{\Gamma}$ bounded uniformly over $\omega$, this number $K_1$ is also independent of $\omega\in[\omega_-,\omega_+]$.

We begin by bounding $\Hnorm{1}{\Omega}{\breve u_0}$ by a multiple of
$\Lnorm{p}{\Omega}{\breve\tau} + \Lnorm{p}{\Omega}{\breve\eps}$.  To do this, we must estimate the right-hand-side of \eqref{approx} in $\Honek^*$.
Since $u$ satisfies the scattering problem \eqref{weak1}, $u$ satisfies the differential equation
\begin{equation}
  \Div \tau\nabla u + \omega^2\eps u = 0
  \quad \text{in } \RR^3.
\end{equation}
Applying Theorem~\ref{thmMeyers} to this equation yields
\begin{equation}
\Lnorm{2}{\Omega}{\omega^2\breve\eps u}
   \leq \omega^2\Lnorm{p}{\Omega}{\breve\eps} \Lnorm{s}{\Omega}{u}
   < C_2 \omega^2 \Lnorm{p}{\Omega}{\breve\eps} \Hnorm{1}{\Omega}{u}
   < C_2 K_1 \omega_+^2 \Lnorm{p}{\Omega}{\breve\eps},
\end{equation}
\begin{multline}
\Lnorm{2}{\Omega}{\breve\tau\nabla u}
  \leq \Lnorm{p}{\Omega}{\breve\tau} \Lnorm{s}{\Omega}{\nabla u} \\
  < C_1\Lnorm{p}{\Omega}{\breve\tau} \left[ \Lnorm{2}{B_1}{u} 
                          + \Lnorm{2}{B_1}{\omega^2\eps u} \right]
  < C_1 K_1 (1+\eps_+^0\omega_+^2) \Lnorm{p}{\Omega}{\breve\tau}.
\end{multline}
From these estimates and Theorem~\ref{thmfieldbound}, we infer that $\breve u_0$ is the unique function in $\Honek$ that satisfies equation \eqref{perturbed} for all $v\in\Honek$ and that there is a constant $K_2$, independent of $\eps$, $\breve\eps$, $\tau$, $\breve\tau$, and $\omega$ such that
\begin{equation}\label{bounded}
  \Hnorm{1}{\Omega}{\breve u_0} < K_2 \left[ \Lnorm{p}{\Omega}{\breve\tau} + \Lnorm{p}{\Omega}{\breve\eps} \right].
\end{equation}
Because of this, the linear functional $(\eps,\tau)\mapsto \breve u_0$ is uniformly bounded from $\left( L^p(\Omega) \right)^2$ to $\Honek$, proving the first part of the Theorem.

An analogous argument can be applied to the system \eqref{perturbed} for $\breve u$
and the corresponding differential equation for $u+\breve u$, which appears on the right-hand side of \eqref{perturbed},
\begin{equation}
  \Div (\tau + \breve\tau)\nabla(u + \breve u) + \omega^2(\eps+\breve\eps)(u + \breve u) = 0
  \quad \text{in } \RR^3.
\end{equation}
This results in the inequality \eqref{bounded}, with $\breve u$ in place of $\breve u_0$, which, together with Lemma \ref{lemmaextension}, yields (reusing the constant $K_2$),
\begin{equation}\label{bounded0}
  \Hnorm{1}{B_2}{\breve u} < K_2 \left[ \Lnorm{p}{\Omega}{\breve\tau} + \Lnorm{p}{\Omega}{\breve\eps} \right].
\end{equation}

Next, we bound the ``second-order" part of the right-hand-side of \eqref{perturbed}, namely $\omega^2\breve\eps\breve u$ and $\nabla\cdot\breve\tau\nabla\breve u$ in $\Honek^*$ by a multiple of $\left( \Lnorm{p}{\Omega}{\breve\tau} + \Lnorm{p}{\Omega}{\breve\eps} \right)^2$.  For this, we apply Theorem~\ref{thmMeyers} to the differential equation
\begin{equation}\label{perturbed2}
\Div \tau\nabla\breve u + \eps\omega^2\breve u = - \Div\breve\tau\nabla(u+\breve u) - \breve\eps\omega^2(u+\breve u) \quad \text{in } B_2.
\end{equation}
\begin{multline}\label{twohatuhat}
\Lnorm{2}{\Omega}{\omega^2\breve\eps\breve u}
  \leq \omega^2 \Lnorm{p}{\Omega}{\breve\eps} \Lnorm{s}{\Omega}{\breve u} \\
  < C_2 \omega^2 \Lnorm{p}{\Omega}{\breve\eps} \Hnorm{1}{\Omega}{\breve u}
  < C_2 K_2 \omega_+^2 \Lnorm{p}{\Omega}{\breve\eps} 
          \left[ \Lnorm{p}{\Omega}{\breve\tau} + \Lnorm{p}{\Omega}{\breve\eps} \right];
\end{multline}
\begin{multline}\label{uno}
\Lnorm{2}{\Omega}{\breve\tau\nabla\breve u}
  \leq \Lnorm{p}{\Omega}{\breve\tau} \Lnorm{s}{\Omega}{\nabla\breve u} \\
  <  C_1 \Lnorm{p}{\Omega}{\breve\tau} \left[
    \Lnorm{2}{B_1}{\breve u} 
         + \omega^2 \Lnorm{2}{B_1}{\eps\breve u} + \omega^2\Lnorm{2}{B_1}{\breve\eps(u+\breve u)}
         + \Lnorm{s}{B_1}{\breve\tau\nabla(u+\breve u)}
    \right].
\end{multline}
The first two terms in the right-hand-side of this last estimate are in turn estimated by
\begin{equation}\label{dos}
  \Lnorm{2}{B_1}{\breve u} + \omega^2 \Lnorm{2}{B_1}{\eps\breve u}
    <  K_2 (1+\omega_+^2\eps_+^0)
      \left[ \Lnorm{p}{\Omega}{\breve\tau} + \Lnorm{p}{\Omega}{\breve\eps} \right].
\end{equation}
The last two terms are estimated by
\begin{multline}\label{tres}
  \omega^2\Lnorm{2}{B_1}{\breve\eps(u+\breve u)}
     \leq \omega^2 \Lnorm{p}{B_1}{\breve\eps} \Lnorm{s}{B_1}{u+\breve u} \\
     < C_2 \omega^2 \Lnorm{p}{B_1}{\breve\eps} \Hnorm{1}{B_1}{u+\breve u}
     < C_2 K_1 \omega_+^2 \Lnorm{p}{B_1}{\breve\eps},
\end{multline}
\begin{multline}\label{cuatro}
  \Lnorm{s}{B_1}{\breve\tau\nabla(u+\breve u)}
  \leq \Lnorm{p}{B_1}{\breve\tau} \Lnorm{q}{B_1}{\nabla(u+\breve u)} \\
  \leq C_1 \Lnorm{p}{B_1}{\breve\tau}
      \left[ \Lnorm{2}{B_2}{u+\breve u} + \Lnorm{2}{B_2}{\omega^2(\eps+\breve\eps)(u+\breve u)}
      \right]
  \leq C_1 K_1 (1+\eps_+^0\omega_+^2) \Lnorm{p}{B_1}{\breve\tau}.
\end{multline}
To simplify notation in the rest of the proof, the symbol $C$ will denote different constants.
By putting the estimates (\ref{uno}, \ref{dos}, \ref{tres}, \ref{cuatro}) together and using the periodicity of $\breve\eps$ and $\breve\tau$ to get
$\left[ \Lnorm{p}{B_1}{\breve\tau} + \Lnorm{p}{B_1}{\breve\eps} \right]
     < C \left[ \Lnorm{p}{\Omega}{\breve\tau} + \Lnorm{p}{\Omega}{\breve\eps} \right]$,
we obtain
\begin{equation}\label{onehatgraduhat}
  \Lnorm{2}{\Omega}{\breve\tau\nabla\breve u}
  < C \Lnorm{p}{\Omega}{\breve\tau} \left[ \Lnorm{p}{\Omega}{\breve\tau} + \Lnorm{p}{\Omega}{\breve\eps} \right].
\end{equation}

Considering $L^2$ functions and their divergences as elements of $H^1(\Omega)^*$, we conclude from \eqref{twohatuhat} and \eqref{onehatgraduhat} that
\begin{equation}\label{H1star}
  \left\| \Div\breve\tau\nabla\breve u + \omega^2\breve\eps\breve u \right\|_{H^1(\Omega)^*}
  <  C \left( \Lnorm{p}{\Omega}{\breve\tau} + \Lnorm{p}{\Omega}{\breve\eps} \right)^2.
\end{equation}
Equations \eqref{perturbed} and \eqref{approx} give the system for $\breve u - \breve u_0$,
\begin{multline}
\int_\Omega
\left( \tau \, \nabla (\breve u-\breve u_0) \cdot \nabla \bar v - \omega^2\, \eps\, (\breve u-\breve u_0) \bar v  \right)
+ \tau \int_{\Gamma} T^\omega (\breve u-\breve u_0) \bar v  \\
\,=\, 
- \int_\Omega
\left( \breve \tau \, \nabla \breve u \cdot \nabla \bar v - \omega^2\, \breve\eps\, \breve u \bar v  \right)
\quad \text{for all } v\in\Honek.
\end{multline}
This, together with Theorem~\ref{thmfieldbound} and equation \eqref{H1star}, gives us the desired result
\begin{equation}
  \Hnorm{1}{\Omega}{\breve u-\breve u_0}
      < C \left( \Lnorm{p}{\Omega}{\breve\tau} + \Lnorm{p}{\Omega}{\breve\eps} \right)^2,
\end{equation}
in which the constant is independent of $\eps$, $\breve\eps$, $\tau$, $\breve\tau$, and $\omega$, subject to the conditions in the Theorem.

Finally, we prove the Lipschitz continuity of the derivative with respect to $\eps$ and $\tau$ in the $L^p$ norm.  To do this, we let $\beps$ and $\btau$ be fixed as directions of differentiation and perturb the functions $\eps$ and $\tau$, at which we differentiate, by functions $\ceps$ and $\ctau$, remembering to require that $\eps$, $\tau$, $\eps+\beps$, $\tau+\btau$, $\eps+\ceps$, $\tau+\ctau$, $\eps+\ceps+\beps$, and $\tau+\ctau+\btau$ satisfy the bounds in the hypotheses of the Theorem.  This results in perturbations of the fields $u$ and $\buz$,
\begin{eqnarray}
  (\eps,\tau) &\mapsto& (\eps+\ceps,\tau+\ctau), \\
  u &\mapsto& u + \cu, \\
  \buz &\mapsto& \buz + \cbuz.
\end{eqnarray}
By subtracting \eqref{approx} as it is written from \eqref{approx} with these substitutions, we obtain
\begin{multline}\label{hihi}
  \int_\Omega ((\tau+\ctau)\nabla\cbuz\cdot\nabla\bar v - \omega(\eps+\ceps)\cbuz\bar v) +
  \tau_0 \int_\Gamma T^\omega\cbuz\bar v \\
  = -\int_\Omega (\btau\nabla\cu\cdot\nabla\bar v - \omega^2\beps\cu\bar v)
  - \int_\Omega (\ctau\nabla\buz\cdot\nabla\bar v - \omega^2\ceps\buz\bar v).
\end{multline}
Steps analogous to those between equations \eqref{twohatuhat} and \eqref{H1star} give us the estimates
\begin{eqnarray}
  \left\| \Div\btau\nabla\cu + \omega^2\beps\cu \right\|_{H^1(\Omega)^*}
  &<&  C \left( \Lnorm{p}{\Omega}{\ctau} + \Lnorm{p}{\Omega}{\ceps} \right)
                \left( \Lnorm{p}{\Omega}{\btau} + \Lnorm{p}{\Omega}{\beps} \right), \\
  \left\| \Div\ctau\nabla\buz + \omega^2\ceps\buz \right\|_{H^1(\Omega)^*}
  &<&  C \left( \Lnorm{p}{\Omega}{\ctau} + \Lnorm{p}{\Omega}{\ceps} \right)
                \left( \Lnorm{p}{\Omega}{\btau} + \Lnorm{p}{\Omega}{\beps} \right),
\end{eqnarray}
in which the constant is independent of the functions.  An application of Theorem~\ref{thmfieldbound} to equation \eqref{hihi} gives a uniform bound
\begin{equation}
  \Hnorm{1}{\Omega}{\cbuz} < C \left( \Lnorm{p}{\Omega}{\ctau} + \Lnorm{p}{\Omega}{\ceps} \right)
                              \left( \Lnorm{p}{\Omega}{\btau} + \Lnorm{p}{\Omega}{\beps} \right).
\end{equation}
This means that the linear functional $(\beps,\btau)\mapsto(\buz+\cbuz)$, which is the derivative of the total field with respect to the material parameters at $(\eps+\ceps,\tau+\ctau)$, is bounded by
$C \left( \Lnorm{p}{\Omega}{\ctau} + \Lnorm{p}{\Omega}{\ceps} \right)$.
We conclude that the derivative, defined on admissible functions $\beps$ and $\btau$, is Lipschitz continuous with respect to the $L^p$ norm over functions $\eps$ and $\tau$ that satisfy the hypotheses of the Theorem.
\end{proof}

\subsection{Transmitted energy}\label{subsectrans} 

We take $b^\text{inc}=0$ so that there is a source field incident upon the slab only from the left.  The energy transmitted to the right-hand side of the slab is given by
\begin{equation}\label{trans}
  \trans = \Im\!\intgp \tau_0\,\bar u\,\partial_nu = -\Im \intgp \tau_0 (Tu)\bar u.
\end{equation}
Let $\breve u$ be a arbitrary perturbation of $u$, and $\breve\trans$ the corresponding perturbation of $\trans$,
\begin{equation}\label{transperturbed}
  \trans + \breve\trans = -\Im \intgp\tau_0 (T(u+\breve u))(\bar u+\bar{\breve u}).
\end{equation}
From the equations for $\trans$ and $\trans + \breve\trans$, we obtain
\begin{equation}
  \breve\trans = -\Im\intgp\tau_0 ((T\breve u) \bar u + (Tu)\bar{\breve u} + (T\breve u)\bar{\breve u}).
\end{equation}
Denote the linear part of $\breve\trans$ by $\breve\trans_0$\,,
\begin{equation}\label{brevetrans0}
  \breve\trans_0(u,\breve u) = -\Im\intgp\tau_0 ((T\breve u) \bar u + (Tu)\bar{\breve u})
\end{equation}
Because of the trace theorem and the boundedness of $T$, we have
\begin{eqnarray}
  && \trans \leq C\Hnorm{1}{\Omega}{u}^2\,, \\
  && |\breve\trans| \leq C\Hnorm{1}{\Omega}{u}\Hnorm{1}{\Omega}{\breve u}\,, \\
  && |\breve\trans - \breve\trans_0| \leq C\Hnorm{1}{\Omega}{\breve u}^2\,.
\end{eqnarray}
This demonstrates that the map $\Honek\to\RR :: u\mapsto\trans$ is bounded and differentiable and that the bounded linear map
\begin{equation}
  \Honek \to \RR :: \breve u \mapsto \breve\trans_0
\end{equation}
defined through \eqref{brevetrans0} is the derivative of $u\mapsto\trans$ at $u$.

Because the adjoint problem for the transmission described below is a scattering problem with wavevector $-\kappa$, we will need to exhibit explicitly the dependence of $\eta_m$, $\zprop$, and $T$ on the wavevector, which we have suppressed until now.  From the definition \eqref{eta} of $\eta_m$ and $\zprop$, we obtain
\begin{equation}
  \etak_{\text{-}m} = \etamk_m\,, \quad \zpropmk = -\zpropk\,.
\end{equation}
Using this and the relation $\hat{\bar u}^\kappa_m = \bar{\hat u}^\kappa_{-m}$ for the $\kappa$-Fourier coefficients \eqref{Fourierbdy} of a function $u\in\Honek$ restricted to $\Gamma$ by the trace map, one can derive the relations
\begin{equation}\label{Tidentity}
  \overline{\Tk^*u} = \Tmk\bar u\,, \quad \overline{\Tk u} = \Tmk^*\bar u\,,
\end{equation}
and thence the equivalent expression for the differential of the transmitted energy
\begin{equation}\label{approxtrans}
  \breve\trans_0 = \breve\trans_0(u,\breve u) = -\Im \intgp \tau_0 ((\Tmk-\Tmk^*)\bar u) \breve u\,.
\end{equation}

Essentially following \cite{LiptonShipmanVenakides2003}, we will demonstrate that the solution $\ua$ to an adjoint scattering problem represents the transmission functional.  We take as incident field a $-\kappa$-pseudoperiodic left-traveling wave $\uainc$, incident upon the slab from the right, that is obtained by sending the transmitted propagating harmonics of $u$ back toward the slab.  This is done by conjugating $u$ and retaining only the propagating harmonics:
\begin{eqnarray}
  && u = \sum_{m\in\ZZ^2} b_m e^{i\etak_mx_3}e^{i(m+\kappa)x'} \quad (x_3\geq z_+), \\
  && \uainc = \text{propagating part of $\bar u$} = \sum_{m\in\zpropmk} \bar b_{\text{-}m} e^{-i\etamk_mx_3}e^{i(m-\kappa)x'}.
\end{eqnarray}
Let $\ua\in\Honemk$ be the solution to Problem \ref{problemscattering2} with Bloch wavevector $-\kappa$ and incident field $\uainc$; thus $\ua$ satisfies
\begin{equation}
  \int_\Omega (\tau \nabla\ua\cdot\nabla\bar w - \omega^2\eps\ua\bar w)
  + \tau_0\!\int_\Gamma(\Tmk\ua)\bar w = \tau_0\!\intgp(\partial_n + \Tmk)\uainc\bar w
  \quad \text{for all } w\in\Honemk.
\end{equation}
Observe that, since $\uainc$ has no right-traveling part,
$(\partial_n + \Tmk)\uainc{|_{}}_{\Gamma_+} = (\Tmk-\Tmk^*)\uainc{|_{}}_{\Gamma_+}$
By making the identification $\breve u_0 = \bar w$, and using the first identity in \eqref{Tidentity} and the definition of $\uainc$ together with the fact that $(\Tmk-\Tmk^*)$ vanishes on the linear and exponential Fourier harmonics, we obtain
\begin{equation}
  \int_\Omega (\tau\nabla\breve u_0\cdot\nabla\ua - \omega^2\eps\breve u_0\ua)
  + \tau_0\int_\Gamma(\Tk\breve u_0)\ua
  = \tau_0 \intgp ((\Tmk - \Tmk^*)\bar u)\breve u_0.
\end{equation}
Using the identification $\ua = \bar v$ in equation \eqref{approx} for the derivative of $(\tau,\eps)\mapsto u$ and equation \eqref{approxtrans}, with $\breve u_0$ in place of $\breve u$, we obtain an expression for the derivative of the composite operation $(\tau,\eps) \mapsto u \mapsto \trans$,
\begin{equation}\label{trans0}
  \breve\trans_0 = \breve\trans_0(u,\breve u_0) = \Im \int_\Omega \left( \breve \tau \, \nabla u \cdot \nabla \ua - \omega^2\, \breve\eps\, u \ua  \right).
\end{equation}

The derivative of $(\tau,\eps) \mapsto u \mapsto \trans$ is Lipschitz continuous, as we now demonstrate.  Consider the derivative at two different pairs $(\eps,\tau)=(\eps^{1,2},\tau^{1,2})$,
\begin{eqnarray}
  (\beps,\btau) \mapsto \buz^1 \mapsto \breve\trans(u^1,\buz^1) & \text{at} & (\eps^1,\tau^1), \\
  (\beps,\btau) \mapsto \buz^2 \mapsto \breve\trans(u^2,\buz^2) & \text{at} & (\eps^2,\tau^2).
\end{eqnarray}
By using equation \eqref{approxtrans}, we obtain
\begin{equation}
 \breve\trans(u^1,\buz^1) - \breve\trans(u^2,\buz^2) =
  -\Im \int_{\Gamma_+}\tau_0\left[(\Tmk-\Tmk^*)(\bar u^1-\bar u^2)\buz^1 + (\Tmk-\Tmk^*)\bar u^2(\buz^1-\buz^2)\right].
\end{equation}
From of the boundedness of $(\Tmk-\Tmk^*)$ and estimates (\ref{bounded0},\ref{bounded},\ref{uniformbound},\ref{Lipschitz}), we obtain the estimate
\begin{equation}
  | \breve\trans(u^1,\buz^1) - \breve\trans(u^2,\buz^2) | \leq
  C \left( \Lnorm{p}{\Omega}{\tau^1\!-\!\tau^2} + \Lnorm{p}{\Omega}{\eps^1\!-\!\eps^2} \right) \left( \Lnorm{p}{\Omega}{\btau} + \Lnorm{p}{\Omega}{\beps} \right),
\end{equation}
which demonstrates the Lipschitz continuity.

The results of this section are summarized in the following theorem.  An analogous theorem can be established for the variational derivative of $b_m$, given by \eqref{b0}.
\begin{theorem}[Transmission sensitivity to $L^p$ perturbations]\label{thmtrans}
  Let the stipulations in Theorem~\ref{thmmain} hold, as well as the equations (\ref{trans},\ref{transperturbed}) defining the transmitted energy $\trans$ and its perturbation $\breve\trans$.  Then the linear operator
\begin{equation}\label{transderiv}
  (L^q(\Omega))^2 \to \RR :: (\breve\tau,\breve\eps) \mapsto \breve\trans_0
\end{equation}
defined through \eqref{trans0} is bounded and
\begin{equation}
  |\breve\trans - \breve\trans_0|
  \leq C \left( \Lnorm{p}{\Omega}{\breve\tau} + \Lnorm{p}{\Omega}{\breve\eps} \right)^2.
\end{equation}
Moreover, the derivative operator \eqref{transderiv} is Lipschitz continuous, restricted to functions $(\eps,\tau)$ that satisfy the hypotheses of Theorem~\ref{thmmain}.
\end{theorem}

\bigskip
\bigskip

\noindent
{\bf \large Acknowledgement.} \ 
This research was supported by the National Science Foundation under grant DMS-0807325.


\bibliography{Perturbation}

\begin{thebibliography}{10}

\bibitem{BaoBonnetier2001}
G.~Bao and E.~Bonnetier.
\newblock Optimal design of periodic diffractive structures.
\newblock {\em Appl. Math. Optim.}, 43:103--116, 2001.

\bibitem{Bao1998}
Gang Bao.
\newblock On the relation between the coefficients and solutions for a
  diffraction problem.
\newblock {\em Inverse Problems}, 14:787--798, 1998.

\bibitem{BendicksoDowlingScalora1996}
Jon~M. Bendickson, Jonathan~P. Dowling, and Michael Scalora.
\newblock Analytic expressions for the electromagnetic mode density in finite,
  one-dimensional, photonic band-gap structures.
\newblock {\em Physical Review E}, 53(4):4107--4121, 1996.

\bibitem{Bonnet-BeStarling1994}
Anne-Sophie Bonnet-Bendhia and Felipe Starling.
\newblock Guided waves by electromagnetic gratings and nonuniqueness examples
  for the diffraction problem.
\newblock {\em Math. Methods Appl. Sci.}, 17(5):305--338, 1994.

\bibitem{ColtonKress1998}
David Colton and Rainer Kress.
\newblock {\em Inverse Acoustic and Electromagnetic Scattering Theory},
  volume~93 of {\em Applied Mathematical Sciences}.
\newblock Springer-Verlag, 2 edition, 1998.

\bibitem{CostabelStephan1985}
Martin Costabel and Ernst Stephan.
\newblock A direct boundary integral equation method for transmission problems.
\newblock {\em J. Math. Anal. Appl.}, 106(2):367--413, 1985.

\bibitem{Dobson1999}
D.~C. Dobson.
\newblock Optimal shape design of blazed diffraction gratings.
\newblock {\em Appl. Math. Optim.}, 40:61--78, 1999.

\bibitem{ElschnerSchmidt1998}
J.~Elschner and G.~Schmidt.
\newblock Diffraction in periodic structures and optimal design of binary
  gratings. part i: direct problems and gradient formulas.
\newblock {\em Math. Meth. Appl. Sci.}, 21(14):1297--1342, 1998.

\bibitem{ElschnerSchmidt2001a}
J.~Elschner and G.~Schmidt.
\newblock Diffraction in periodic structures and optimal design of binary
  gratings ii: Gradient formulas for tm polarization.
\newblock {\em Prob. Meth. in Math. Phys: the Siegfried Pr\"ossdorf Memorial
  Volume (Oper. Theory: Adv. Appl.)}, 121:89--108, 2001.

\bibitem{ElschnerSchmidt2003}
J.~Elschner and G.~Schmidt.
\newblock Conical diffraction by periodic structures: Variation of interfaces
  and gradient formulas.
\newblock {\em Math. Nachr.}, 252:24--42, 2003.

\bibitem{ElschnerSchmidt2001}
Johannes Elschner and Gunther Schmidt.
\newblock Inverse scattering for periodic structures: stability of polygonal
  interfaces.
\newblock {\em Inverse Problems}, 17:1817--1829, 2001.

\bibitem{GilbargTrudinger1998}
David Gilbarg and Neil~S. Trudinger.
\newblock {\em Elliptic Partial Differential Equations of Second Order}.
\newblock Springer-Verlag, 1998.

\bibitem{Kirsch1994}
Andreas Kirsch.
\newblock Uniqueness theorems in inverse scattering theory for periodic
  structures.
\newblock {\em Inverse Problems}, 10:145--152, 1994.

\bibitem{Kuchment1993}
Peter Kuchment.
\newblock {\em Floquet Theory for Partial Differential Equations}.
\newblock Birkh{\"a}user Verlag AG, 1993.

\bibitem{LenoirVulliermeHazard1992}
M.~Lenoir, M.~Vullierme-Ledard, and C.~Hazard.
\newblock Variational formulations for the determination of resonant states in
  scattering problems.
\newblock {\em SIAM J. Math. Anal.}, 23(3):579--608, 1992.

\bibitem{LiptonShipmanVenakides2003}
Robert~P. Lipton, Stephen~P. Shipman, and Stephanos Venakides.
\newblock Optimization of resonances in photonic crystal slabs.
\newblock In Philippe Lalanne, editor, {\em Physics, Theory, and Applications
  of Periodic Structures in Optics II}, volume 5184 of {\em Proceedings
  Series}, pages 168--177. SPIE--The International Society for Optical
  Engineering, 2003.

\bibitem{Meyers1975}
Norman Meyers.
\newblock An {$L^p$} estimate for the gradient of solutions of second order
  elliptic divergence equations.
\newblock {\em Duke Mathematical Journal}, 42:121--136, 1975.

\bibitem{MoussaWangTuttle2007}
R.~Moussa, B.~Wang, G.~Tuttle, Th. Koschny, and C.~M. Soukoulis.
\newblock Effect of beaming and enhanced transmission in photonic crystals.
\newblock {\em Phys. Rev. B}, 76:235417--1--8, 2007.

\bibitem{MuratTartar1978}
Fran{\c c}ois Murat and Luc Tartar.
\newblock {\em H-Convergence}, volume 31: Topics in the Mathematical Modelling
  of Composite Materials of {\em Progress in Nonlinear Differential Equations
  and Their Applications}.
\newblock Birkh{\"a}user Verlag AG, 1978.

\bibitem{Neviere}
M.~Nevi{\`e}re.
\newblock {\em The homogeneous problem, in Electromagnetic Theory of Gratings},
  chapter~5, pages 123--157 (Ed. R. Petit).
\newblock Springer-Verlag, Berlin, 1980.

\bibitem{Pironneau1983}
Olivier Pironneau.
\newblock {\em Optimal Shape Design for Elliptic Systems}, volume xiii of {\em
  Series in Computational Physics}.
\newblock Springer-Verlag, 1983.

\bibitem{RamdaniShipman2008}
Karim Ramdani and Stephen~P. Shipman.
\newblock Transmission through a thick periodic slab.
\newblock {\em Math. Mod. Meth. Appl. S.}, 18(4):543--572, 2008.

\bibitem{ReedSimon1980d}
Michael Reed and Barry Simon.
\newblock {\em Methods of Mathematical Physics: Analysis of Operators},
  volume~IV.
\newblock Academic Press, 1980.

\bibitem{ShipmanVenakides2005}
Stephen~P. Shipman and Stephanos Venakides.
\newblock Resonant transmission near non-robust periodic slab modes.
\newblock {\em Physical Review E}, 71(1):026611--1--10, 2005.

\bibitem{ShipmanVolkov2007}
Stephen~P. Shipman and Darko Volkov.
\newblock Guided modes in periodic slabs: existence and nonexistence.
\newblock {\em SIAM J. Appl. Math.}, 67(3):687--713, 2007.

\bibitem{Wilcox1984}
Calvin~H. Wilcox.
\newblock {\em Scattering Theory for Diffraction Gratings}, volume~46 of {\em
  Applied Mathematical Sciences}.
\newblock Springer-Verlag, 1984.

\end{thebibliography}

\end{document}